\documentclass{article}
\usepackage{arxiv}
\usepackage{amsmath,amsfonts,amssymb,amsthm}
\usepackage{mathHeader}

\usepackage[utf8]{inputenc} % allow utf-8 input
\usepackage[T1]{fontenc}    % use 8-bit T1 fonts
\usepackage{csquotes}
\usepackage{hyperref}       % hyperlinks
\usepackage{url}            % simple URL typesetting
\usepackage{booktabs}       % professional-quality tables
\usepackage{amsfonts}       % blackboard math symbols
\usepackage{nicefrac}       % compact symbols for 1/2, etc.
\usepackage[protrusion=true,expansion=false]{microtype}      % microtypography
\usepackage{hyphenat}

\usepackage{commath} % for \dif command
\usepackage{colortbl}
\usepackage{pgfplots}
\pgfplotsset{compat=newest}

\definecolor{cell}{rgb}{0.9,0.9,0.9}

\graphicspath{{../}}
\makeatletter
\def\input@path{{../}}
\makeatother

\title{Considering Slow Manifold Based Model Reduction for Multiscale Chemical Optimal Control Problems}
\author{
	Marcus~Heitel\\ %\thanks{Use footnote for providing further information about author (webpage, alternative 	address)---\emph{not} for acknowledging funding agencies.} \\
	Institute for Numerical Mathematics\\
	Ulm University\\
	Ulm, Germany
	\And
	Robin~Verschueren \\
	Department of Microsystems Engineering (IMTEK)\\
	University of Freiburg\\
	Freiburg, Germany
	\And
	Dirk~Lebiedz \\
	Institute for Numerical Mathematics\\
	Ulm University\\
	Ulm, Germany
	\And
	Moritz~Diehl \\
	Department of Microsystems Engineering (IMTEK)\\
	University of Freiburg\\
	Freiburg, Germany
}

\begin{document}
\maketitle

\begin{abstract}
	Finite-dimensional dissipative dynamical systems with multiple time-scales are obtained when modeling chemical reaction kinetics with ordinary differential equations. Such stiff systems are computationally hard to solve and therefore, optimal control problems which contain chemical kinetic models as infinitesimal constraints are even more difficult to handle. Model reduction might offer an approach to improve numerical efficiency as well as avoid stiffness of such models. We show in this paper for benchmark problems how attracting manifold computation methods could be exploited to solve optimal control problems more efficiently while having in mind the ambitious long-term goal to apply them to real-time control problems in chemical kinetics.
\end{abstract}

% keywords can be removed

\section{Introduction}
\label{sec:intro}
Modeling chemical reaction mechanisms with ordinary differential equations (ODEs) yields highly-nonlinear systems with multiple 
time-scale dynamics. This makes them computationally hard to control online, e.g. with nonlinear model predictive control (NMPC). 
Thus, model reduction methods can be considered to tackle such problems. 

A class of model reductions that can used in the context of multiple time-scale dynamics is the class of manifold-based model reduction methods (for an overview see e.g. Lebiedz and Unger\cite{Lebiedz2016}). 
Since the long time behavior of such systems is dominated by the slow components which are often stable, the system usually approaches low\hyp dimensional manifolds in the phase space divided into slow and fast components, i.e. the slow invariant manifold is attracting.

It is still an open issue how these manifolds can be used to simplify and accelerate the computations of optimal control problems. First investigations
towards this goal have been performed by Rehberg \cite{Rehberg2013, Lebiedz2013b}. % {\color{red} auch arXiv zitieren}.
In this paper, we show how (point-wise) online computation of such manifolds can be integrated into the multiple shooting based solution of optimal control problems in a reasonable and useful way. 
Two low\hyp dimensional benchmark problems demonstrate the efficiency of this method.

%	In this paper, we show how (point-wise) online computation of such manifolds can be integrated in a reasonable and useful way into the solution of optimal control
%	problems and compare these approaches with already existing direct numerical optimal control approaches (without manifold computation) such as the shooting methods 
%	and collocation.

\subsection{Classification of Model Reduction Methods}
Modern model order reduction methods commonly in use for optimal control, in particular in large-scale partial differential equation (PDE) applications, can roughly be divided into at least three categories: Krylov based methods, singular value decomposition (SVD) based methods and trajectory based methods. In particular, the latter is a nonlinear method, whereas the former two have originally been developed for large linear systems.
We give a brief description of them and list some of their main characteristics, advantages as well as disadvantages.

\textbf{Krylov Based Methods} project systems onto Krylov subspaces. A widely used technique of this category 
is the moment matching method \cite{Feldmann1995} which interpolates 
derivatives (moments) of the transfer function for linear, time-independent systems at certain points. Arnoldi- 
or Lanczos processes are used for the projection onto the Krylov subspaces and can be implemented efficiently. But there exists no global
error bound.

In \textbf{SVD Based Methods} high-dimensional, high-rank matrices $M$ are approximated by matrices $M_{k}$ of lower rank $k$ using singular value decomposition. 
Furthermore, these rank $k$ approximations are optimal in unitarily invariant norms and for the Euclidean norm. The replacement of the original matrix and its approximation is determined by the $(k+1)$-th largest singular value. 
A very popular SVD based method is balanced truncation \cite{Moore1981}. It balances linear, time-independent systems such that their controllability
and observability Gramians become diagonal and equal. Negligible states are truncated such that asymptotic stability can be preserved. A priori error bounds are computable. 

For nonlinear and time-dependent systems the probably most successful method
is proper orthogonal decomposition (POD). The key idea of POD is to compute samples of the trajectory - so called \emph{snapshots} - 
and use SVD in order to construct a basis of a subspace of the snapshot space. 
However, as described in Chaturantabut and Sorensen\cite{Chaturantabut2010}, the computational complexity of the reduced model is almost as high as for 
the original model, if the system is very nonlinear. To overcome this difficulty, for example Empirical Interpolation Methods \cite{Chaturantabut2010},\cite{Drohmann2012} are used.
For a more detailed overview of the SVD based methods and the Krylov based methods see Antoulas and Sorensen\cite{Antoulas2001} or Benner, Sachs, and Volkwein\cite{Benner2014}.

\textbf{Trajectory Based Methods} are totally different in nature, they are fully nonlinear, aimed at a replacement of the full (usually linear) state space by a (nonlinear) manifold of lower dimension. They are based on time scale separation which arises for 
example in models of chemical kinetics. The state space is not projected onto a linear or affine subspace but onto a lower dimensional 
manifold. So the nonlinearity of the problem is somehow encoded in the manifold. 
Maas and Pope \cite{Maas1992} introduced the intrinsic low\hyp dimensional manifolds (ILDM) method which uses the structure 
of a spectral decomposition of the Jacobian of the vector field. Decoupling slow and fast submatrices yields 
a good approximation of the slow invariant manifold (SIM). The reduced model also has the benefit of less stiff dynamics.
A drawback of this method and manifold-based approaches in general is that there are no reliable a priori error bounds.

In this paper, we restrict ourselves to a subclass of trajectory based methods based on some of the authors' previous work. Our model order reduction approach \cite{Lebiedz2004,Reinhardt2008,Lebiedz2010,Lebiedz2010a,Lebiedz2011,Lebiedz2013} is based on either a nonlinear root finding problem or a equation of a variational principle for approximation of trajectories on slow (attracting) invariant manifolds. 

Both POD and trajectory based methods exploit a contraction property of the dynamical systems for the construction subspaces. We choose a trajectory based method for optimal control in the context of chemical reactions because of the high nonlinearity and stiffness that come with chemical dynamics.

The paper is organized as follows. A brief introduction to manifold-based model reduction is given in Section \ref{sec:modelreduction} 
and two successful numerical model reduction methods are presented. A new approach to combine
manifold-based calculation of slow manifolds and optimal control problems is introduced. 
Numerical methods for the solution of optimal control problems can be found in Section \ref{sec:num}, results are in Section \ref{sec:results}, 
and a conclusion follows in Section \ref{sec:conclusion}.

\section{Model Reduction Methods for Multi-Scale Systems} 
\label{sec:modelreduction}
Let a general system be given by
\begin{align}\label{eq:ode}
	\dot{z}(t) = f\left(z(t)\right) \quad \forall t \in[0,T],
\end{align}
where $f \in C^r(\R^n,\R^n)$ for a number $r \in \N \cup \lbrace \infty \rbrace$, i.e. $f$ is a $r$-times continuously differentiable function, that maps from $\R^n$ to $\R^n$.
%{\color{red} reicht das ? $C^1$ ?, $C^\infty$ ?}
In the following we assume that system \eqref{eq:ode} has multiple time-scales, i.e.
the eigenvalues $\lambda_1 \le \hdots \le \lambda_n$ of the Jacobian of $f$ at the equilibrium (a root of $f$) are clustered in the sense that there exists a $1\le k \le n-1$ such that $\lambda_1 \le \hdots \le  \lambda_k \ll \lambda_{k+1} \le \hdots \le \lambda _n$.

In the phase space of such systems often there exists a manifold which is invariant under the flow of the dynamics \eqref{eq:ode} and 
that is characterized by a fast bundling of trajectories onto this manifold. For a given initial value, the corresponding
trajectory will reach the neighborhood of such a manifold fast and then stay close to it. Therefore, this manifold is called slow invariant attractive manifold (SIM).

For simplicity, we only consider systems in explicit singularly perturbed form, i.e.
\begin{subequations}\label{eq:sps}
	\begin{alignat}{2}
	\dot{z_s}(t) &= f_s(z_s(t),z_f(t); \epsilon)\\
	\varepsilon \dot{z_f}(t) &= f_f(z_s(t),z_f(t); \epsilon) \quad \forall t \in[0,T],
	\end{alignat}
\end{subequations}
with $0 < \epsilon \ll 1$, $f_s \in C^r(\R^{n_s}\times \R^{n_f}, \R^{n_s})$ and $f_f \in C^r(\R^{n_s}\times \R^{n_f}, \R^{n_f})$.
Here, $\epsilon$ represents a measure for the time-scale separation. The functions $f_s$ and $f_f$ may depend on $\epsilon$ in a polynomial way.

The analysis of multi-scale systems is investigated well for singularly perturbed systems.
Fenichel \cite{Fenichel1979} showed the existence of the SIM in the form of Theorem \ref{theo:Fenichel}.

\begin{theorem}[Fenichel \cite{Fenichel1979},\cite{Kuehn2015}]\label{theo:Fenichel}
	Let $S_0 = \left\{(z_s,z_f): f_f(z_s,z_f;0)=0\right\}$ be compact, attracting under flow of \eqref{eq:sps} and $f_s,f_f \in C^r$ for a natural number $r<\infty$, then it holds:
	$\exists \; \epsilon_0 > 0 $ such that $\forall\;0<\epsilon\le \epsilon_0,$  there is a mapping $h_\epsilon(\cdot):K \subset \R^{n_s} \ra \R^{n_f}$ with
	$$ S_\epsilon := \left\{(z_s,z_f): z_f = h_\epsilon(z_s), z_s \in K\right\}.$$
	$S_\epsilon$ is attracting and locally invariant under the flow of \eqref{eq:sps}. Further, $h_\epsilon(\cdot)$ has the asymptotic expansion
	$$ h_\epsilon(z_s) = \sum_{i=0}^{r}h_{\epsilon,i}(z_s)\epsilon^i + \mathcal{O}\left(\epsilon^{r+1}\right)$$
	and Hausdorff distance $\mathcal{O}(\epsilon)$ from $S_0$.
\end{theorem}

\begin{remark}
	As the fast variables $z_f$ can be reconstructed out of the slow variables $z_s$, if the manifold mapping 
	$h_\epsilon$ is known, the slow variables are the so called \emph{reaction progress variables (RPVs)}.
	
	In the following, $(z_s(t),z_f(t))$ will often be abbreviated by $z(t)$ and correspondingly 
	$(f_s,f_f)$ will be abbreviated by $f$.
\end{remark}

Theorem \ref{theo:Fenichel} is the theoretical basis and motivation of manifold-based model reduction for singularly perturbed systems since system \eqref{eq:sps} can be reduced to
\begin{alignat}{2}\label{eq:spsReduced}
\dot{z_s}(t) &= f_s\left(z_s(t),h_\epsilon(z_s(t)); \epsilon\right) \quad \forall \, t \in[0,T].
\end{alignat}
This reduced system \eqref{eq:spsReduced} is of smaller dimension than \eqref{eq:sps} and our calculations in Section \ref{sec:results} show that 
it is less stiff because fast time scales have been eliminated. Especially, the second property is numerically important as we will see in Section \ref{sec:results}. 
However, Theorem \ref{theo:Fenichel} does not provide any efficient method for the calculation of $h_{\epsilon}$.
%	(see Deuflhard and Heroth\cite{Deuflhard1996})

\subsection{Methods for the Approximate Calculation of the SIM}
\label{sec:SIM_calc}
In the following we will shortly present numerical methods that are able to calculate points of the SIM approximately and that proved to be successful in slow manifold computation, i.e. approximations 
for the function $h_\epsilon(\cdot)$. These methods 
represent functions $h_{\text{app}}:\R^{n_s} \rightarrow \R^{n_f}$ such that the point $\left(z_s^*,h_{\text{app}}(z_s^*)\right)$ is located 
in a close neighborhood of the SIM. The error $e(t):= \norm{z_f(t)-h_{\text{app}}(z_s(t))}=\norm{h_\epsilon(z_s(t))-h_{\text{app}}(z_s(t))}$ should be 
as small as possible, however in general we have no error estimates. In case of stable systems this usually does not give rise to severe problems in applications since the error decreases exponentially with time while the dynamical system state evolves with the flow.

\vspace*{1em}
{\bfseries Approach Lebiedz/Unger \cite{Lebiedz2016}}
The first approach suggested by Lebiedz and Unger is motivated mainly by the following fact: Among arbitrary trajectories of system \eqref{eq:ode} 
for which the slow components end within the time $t_1-t_0$ in the state $z_s^*$ a corresponding trajectory piece close the SIM is characterized by 
the smallest curvature in time-parametrization (for motivation, theoretical justification and applications see Lebiedz et al. 
\cite{Lebiedz2004,Reinhardt2008,Lebiedz2010,Lebiedz2010a,Lebiedz2011,Lebiedz2013}).
This motivates optimization problem \eqref{eq:Unger} in terms of a boundary value problem (BVP).

\begin{subequations}\label{eq:Unger}
	\begin{align}
		&\underset{z(\cdot)=\left(z_s(\cdot),z_f(\cdot)\right)}{\text{minimize}}	 &  & \norm{\ddot{z}(t_0)}_2^2 & &\\
		& \; \;\text{such that} & \dot z(t) &= f\big(z(t)\big),   & &\text{for}~t \in [t_0,t_1],	\\
		&						 & 	0			& \le c(z(t)), & & \\
&						 &	z_s^*	&= z_s(t_1). & &
	\end{align}
\end{subequations}

Function $c \in C^\infty(\R^n, \R^\ell)$ contains possible additional constraints, e.g. mass conservation conditions for the 
species.
Here, $h_{\text{app}}(z_s^*)$ is defined as the fast components of the $\arg \min$ of \eqref{eq:Unger}. 
No general estimation of the error $e(t)$ is known. The choice of $t_0 = t_1$ is computationally easier but less accurate. We call
that special choice the local equation of problem \eqref{eq:Unger}.
An advantage of the local approach is that the BVP degenerates into an unconstrained nonlinear program
and thus, it can be handled with a root finding problem via the first order conditions.

\vspace*{1em}
{\bfseries Zero Derivative Principle (ZDP) \cite{Gear2005}}
The idea for the Zero Derivative Principle is to eliminate fast dynamical modes by zeroing higher-order derivatives of the vector field. For a given slow components vector $z_s^*$, 
the ZDP identifies a point $z_f^*$ in terms of a root finding problem: 
\begin{align}\label{eq:ZDP}
	\frac{\dif{}^m z_f}{\dif t^m} = 0 \qquad \text{for a given } m \in \N.
\end{align}
For example, if $m=2$ the ZDP equation \eqref{eq:ZDP} can be rewritten to
\begin{align}\label{eq:ZDPm=2}
	\frac{\dif{}^2 z_f}{\dif t^2} = \frac{\dif}{\dif t}f_f(z_s^*,z_f;\epsilon) = \frac{\partial f_f(z_s^*,z_f;\epsilon)}{\partial z_f}\cdot
	f_f(z_s^*,z_f;\epsilon) 
\end{align}	
using the chain rule and is thus very similar to \eqref{eq:Unger} in the local approach. Further, \eqref{eq:ZDPm=2} can be seen as 
a root finding problem of the function 
$$\psi(z_f) = \frac{\partial f_f(z_s^*,z_f;\epsilon)}{\partial z_f}\cdot
f_f(z_s^*,z_f;\epsilon). $$

Zagaris et al. \cite{Zagaris2009} show that the difference between the SIM 
(given by Fenichel's definition in terms of an $\epsilon$-series) and $h_{\text{app}}$ (here: root of $\psi$) given by the ZDP is $e(t) = \mathcal{O}(\epsilon^m)$.

The ZDP was generalized by Beno\^{\i}t, Br{\o}ns, Desroches and Krupa \cite{Benoit2015}. So, the SIM can also be approximated by setting the 
time-derivatives of the slow components to zero. If the derivative of the $i$-th component of $f_s$ w.r.t. $z_s$
for a given point $(z_s^*,z_f^*)$ is regular then the equation \eqref{eq:generalZDP}
\begin{align}\label{eq:generalZDP}
\frac{\dif{}^m z_s}{\dif t^m} = 0 \qquad \text{for a given } m \in \N
\end{align}
gives an approximation of the SIM up to order $\mathcal{O}(\epsilon^{m-1}) = e(t)$.

Since the SIM is intrinsically characterized by the dynamics, it should be possible to calculate the SIM
for an arbitrary choice of the slow variables (thus, the SIM should be represented in coordinate-independent form). 
Therefore, a reasonable numerical method should also be invariant under such a choice 
of \enquote{slow variable parametrization}. The generalized ZDP shows that this is possible in principle. However, the generalized ZDP method is strictly speaking not coordinate-free,
because an explicit choice of slow variables has to be made and the numerical result slightly depends on this choice. 

It should be mentioned that the ZDP method is local, i.e. only in a neighborhood of $(z_s^*,z_f^*,0)$ the ZDP method
$$\frac{\dif{}^m z_f}{\dif t^m} = 0$$ yields an error approximation  $e(t) = \mathcal{O}(\epsilon^m)$. 
A priori it is not clear how large this neighborhood is. So, in practice one has to carefully check and evaluate the numerical solution.

\subsection{Slow Manifold Based Optimal Control}
Consider ODE constrained optimal control problems (OCP) in the following form:

\begin{subequations}\label{eq:origOCP}
	\begin{align}
		&\underset{z_s(\cdot),\,z_f(\cdot),\,u(\cdot)}{\text{minimize}}	 & \multicolumn{2}{c}{$\displaystyle  \int_{0}^{T}L(z_s(t),z_f(t),u(t)) \dif t $ } \\
%				&\underset{z_s(\cdot),\,z_f(\cdot),\,u(\cdot)}{\text{minimize}}	 &  & \int_{0}^{T}L(z_s(t),z_f(t),u(t)) \dif t  & &\\
		& \; \text{such that} &\dot{z}_s(t) &=f_s(z_s(t),z_f(t),u(t);\epsilon),   & &\text{for}~t \in [0, T],	 \label{eq:ode1InOrigOCP}\\
		&						 & \epsilon \dot{z}_f(t) &=  f_f(z_s(t),z_f(t),u(t); \epsilon), & & \text{for}~t \in [0, T], \label{eq:ode2InOrigOCP}\\
		&						 & 	0			& \le c(z(t)), & & \\
		&						 & z_s^{(0)}	&= z_s(0),& & \\
		&						 & z_f^{(0)}	&= z_f(0). & &
	\end{align}
\end{subequations}

where $z_s(t) \in \R^{n_s}, z_f(t) \in \R^{n_f}$ with control $u(t) \in \R^{n_u}$, for $t \in [0, T]$. In the following, we will refer to~\eqref{eq:origOCP} as the original or full problem. We will use it as reference problem for the comparisons of Section~\ref{sec:results}.

In realistic chemical reaction mechanisms this full OCP often cannot be solved in real-time, because either the total number of
optimization variables or the stiffness of ODE~\eqref{eq:sps} prevents numerical algorithms from calculating local 
optima of~\eqref{eq:origOCP} efficiently.

A possibility to reduce its dimension is to insert the manifold mapping $h_\epsilon$ in \eqref{eq:origOCP} directly \cite{Rehberg2013}:

\begin{subequations}\label{eq:reducedOCP}
	\begin{align}
		&\underset{z_s(\cdot),u(\cdot)}{\text{minimize}}	 &  \multicolumn{2}{c}{$ \displaystyle \int_{0}^{T}L\left(z_s(t),h_\epsilon\left(z_s(t),u(t)\right),u(t)\right) \dif t $}\\
		& \text{such that} &\dot{z}_s(t) &= f_s\left(z_s(t),h_\epsilon\left(z_s(t),u(t)\right),u(t);\epsilon \right),   & &\text{for}~t \in [0, T],	\label{eq:odeInRedOCP}\\
		&						 &z_s^{(0)}	&=  z_s(0), & & 
	\end{align}
\end{subequations}

where again $h_\epsilon(z_s)$ can be numerically approximated by $h_{\text{app}}$ resulting from the Lebiedz/Unger or the ZDP approach. 

Note that the controls $u(t)$ are regarded as \enquote{slow} variables, because of several reasons. In chemical reactions the fast time scale is often too fast for precise control of 
input species. In the next section we discretize $u(t)$ and the other states in time and use piecewise constant values to approximate $u(t)$.
This prevents the optimal control from chattering with a huge number of switches and constant values are directly related to slow variables.

The dynamics of \eqref{eq:odeInRedOCP} are less stiff than those of \eqref{eq:origOCP} and 
\eqref{eq:reducedOCP} has just $n_s + n_u$ instead of $n_s + n_f + n_u$ optimization variables.
Numerically however, $h_\epsilon(\cdot)$ in problem \eqref{eq:reducedOCP} will be approximated by a root finding respectively optimization problem, which 
has to be solved in each iteration of any algorithm for optimization problem \eqref{eq:reducedOCP}.

\vspace*{1em}
{\bfseries A New Approach for Efficiently Solving Manifold-Based OCPs}

The new approach we propose is based on the idea of \emph{lifting}, as first presented in Albersmeyer and Diehl\cite{Albersmeyer2010}.
To prevent numerical solvers from calculating $h_\epsilon(z_s,u)$ respectively $h_\text {app}(z_s,u)$ exactly in each iteration 
of \eqref{eq:reducedOCP} the equation $z_f=h_\epsilon(z_s,u)$ is handled as a constraint of the OCP. The numerical solution
has to satisfy this equation, but not other iterates, which cuts computational costs.
In comparison with the original OCP \eqref{eq:origOCP}, the dynamics \eqref{eq:ode1InOrigOCP} and \eqref{eq:ode2InOrigOCP} are replaced 
by the following differential algebraic equation (DAE) \eqref{eq:DAE}:
\begin{subequations}\label{eq:DAE}
	\begin{alignat}{2}
	\dot{z}_s &= f_s(z_s, z_f,u;\epsilon) \\
	0 &= z_f - h_\epsilon(z_s,u).
	\end{alignat}
\end{subequations}
Using the methods of Lebiedz/Unger or ZDP for the approximation of $h_\epsilon$ (both can be written as root finding problems $\psi=0$), 
we get the following OCP 
\eqref{eq:liftedOCP}:

\begin{subequations}\label{eq:liftedOCP}
	\begin{align}
		&\underset{z_s,z_f,u}{\text{minimize}}	 &  & \int_{0}^{T}L(z_s,z_f,u) \dif t & &\\
		& \text{such that} &\dot{z_s} &= f_s(z_s,z_f,u;\epsilon),   & &   \label{eq:dae1InLiftedOCP}\\
		&						 & 	0			&= \psi(z_s, z_f, u):=z_f-h_\text{app}(z_s,u),  \label{eq:dae2InLiftedOCP}& & \\
		&						 & z_s^{(0)}	&= z_s(0),& & \\
		&						 & z_f^{(0)}	&= z_f(0). & &
	\end{align}
\end{subequations}

In the following, we will refer to \eqref{eq:reducedOCP} as the reduced OCP and to \eqref{eq:liftedOCP} as the lifted OCP. 
It has the same number of optimization variables like the full OCP \eqref{eq:origOCP} and its dynamics \eqref{eq:dae1InLiftedOCP} are less stiff such 
that the usage of an explicit integrator is often justified.

\section{Numerical Solution} \label{sec:num}

The choice of numerical solution method for OCPs~\eqref{eq:origOCP}, \eqref{eq:reducedOCP}, \eqref{eq:liftedOCP} is an important decision to make in order to target real-time applications.

We will solve said OCPs with a direct approach. More precisely, we choose a multiple shooting discretization~\cite{Bock1984}, which has been proven to work well and efficiently in practice~\cite{Diehl2002}. We take $N$ shooting intervals, on which we simulate the slow dynamics using an explicit Runge-Kutta integrator of order $4$, where we take $1$ integration step per shooting interval. Note that the explicit integration method is only suited for the reduced and lifted problem formulations. For the full problem, we use an implicit integration scheme, namely a stiffly accurate Radau II-A method. For all methods discussed, we model the controls as piecewise constant over one shooting interval. %We remark that no DAE solver is needed, since we already possess a function $h_\epsilon$ which gives us the fast variables $z_f$.

We noted before that the new efficient method for solving manifold-based OCPs is compatible with either the Lebiedz/Unger approach or with the ZDP. In the latter case, we are free to pick any $m$, such that the $m$-th time derivative should be equal to zero (see Section~\ref{sec:SIM_calc}). We chose $m=2$, as it is a good trade-off between accuracy (recall that the error is $\mathcal{O}(\epsilon^m)$) and computational cost of forming the derivatives.

Carrying out the discretization yields a finite-dimensional nonlinear programming problem (NLP), which we then solve with a Newton-type optimization method. More specifically, we use the interior point solver IPOPT~\cite{Waechter2006}, which we call from within the CasADi~\cite{Andersson2013} framework for dynamic optimization. CasADi also implements the necessary algorithmic differentiation routines for the higher-order derivatives of the functions present in the OCP, such that we are able to use exact Hessians in our solution method.

%	We remark that different choices for discretization of the problem exist. For example, we can employ direct collocation instead of multiple shooting. This  discretization method yields NLPs with even higher dimensions. Nevertheless, we also use direct collocation for comparison, which is only useful if the dynamics do not allow the usage of explicit integrators. Otherwise it yields just an overhead of variables. Therefore, we only solve the full OCP with direct collocation, but not the reduced or lifted OCP. 
We remark that different choices for discretization of the problem exist. The multiple shooting method is chosen, because it widely spread for optimal control problems with PDE and ODE constraints. Its advantages are numerical stability (e.g., in contrast with the single shooting) and the usage of a moderate number of optimization variables. Therefore, it is often used for high-dimensional optimal control problems.
	
For readers that are familiar with CasADi it should be mentioned that in the implementation of the different approaches \texttt{SX} variables are used whenever it is possible and \texttt{MX} variables are only used for implicit integrators. Symbolic \texttt{SX} expressions are typically longer, but faster to evaluate in comparison with \texttt{MX} expressions. However, \texttt{MX} expressions are more flexible with a larger range of applications, e.g., for root-finding problems based on Newton's method. 

\section{Results} \label{sec:results}

For all numerical experiments an Ubuntu 18.04 machine with an Intel\textregistered~Core\texttrademark~i7-2600 (8 cores with 3.40 GHz)
CPU and 16~GB RAM is used. The software version R2018a of MATLAB and version 3.4.4 of CasADi are installed. 

The following examples serve as first benchmark problems and do not have a deeper meaning or specific aim in mind. The objective functions are artificial.
\subsection{Enzyme Example} A simple chemical reaction mechanism with enzymes is given by the Michaelis-Menten-Henri mechanism 
(cf. example 11.2.4 in Kuehn\cite{Kuehn2015})
\begin{align}\label{eq:enzyme}
S + E \rightleftharpoons C \rightarrow P + E,
\end{align}
where $S$ is a substrate, $E$ is an enzyme, $C$ the corresponding substrate-enzyme-complex and $P$ a product.
Simplifying the ODE \cite{Fall2002} which underlies \eqref{eq:enzyme} and introducing an artificial objective function yields
optimal control problem \eqref{eq:enzymeOCP}.

\begin{subequations}\label{eq:enzymeOCP}
	\begin{align}
		&\underset{z(\cdot),u(\cdot)}{\text{minimize}}	 &  & \int_{0}^{5} -50\, z_f(t) + u^2(t) \dif t & &\\
		& \text{such that} &\dot{z_s}(t) &=-z_s(t) +\left(z_s(t) + 0.5\right)z_f(t) + u(t),   & & \text{for}~t \in[0,5], \\
		&						 & \epsilon \dot{z_f}(t) 			&= z_s(t)- \left(z_s(t)+1\right)z_f(t)  ,& & \text{for}~t \in[0,5],\\
		&						 & 1	&= z_s(0),& & \\
	\end{align}
\end{subequations}

where $z(\cdot) = \left(z_s(\cdot),z_f(\cdot)\right) \in \R_+ \times \R_+$ and control $u(t) \in [0,10]$ represents the possibility to add
some substrate (corresponds to variable $z_s$) to the system.
	
The numerical results for this example discretized in $N=40$ subintervals and with $\epsilon = 10^{-6}$ can be seen in 
Figures \ref{fig:enzymeFull},\ref{fig:enzymeZDP},\ref{fig:enzymeDiffAbs}, and \ref{fig:enzymeDiffRel}. 
The solutions of the original system \eqref{eq:origOCP} and
lifted system \eqref{eq:liftedOCP} are almost equal. To be more precise, it holds 
\[
	\max\left\lbrace \lvert \lvert 	z_s^{\text{orig}}-z_s^{\text{app}} \rvert \rvert_{\infty}, 
	\lvert \lvert	z_f^{\text{orig}}-z_f^{\text{app}} \rvert \rvert_{\infty},
	\lvert \lvert	u^{\text{orig}}-u^{\text{app}} \rvert \rvert_{\infty} \right\rbrace \nonumber 
	= \lvert \lvert 	z_s^{\text{orig}}-z_s^{\text{app}} \rvert \rvert_{\infty} \approx 0.047,
\]
which gives a relative error smaller than 1.16 \% (except of the last control $u_N$ that is almost zero). The objective function values produce a relative error
of 0.21 \% (the objective function value is about -187.85). 
For brevity, we omit the numerical solution of the direct collocation method here. The errors between the multiple shooting and the collocation solution are 
comparable to the errors between the multiple shooting of the full OCP and the lifted OCP. Just the difference in the numerical values for $z_s$ is a little bit bigger.
To sum up, there is no big difference in the quality of the solutions.

\begin{figure}[ht]
	\centering
	\bmp{0.4}
%\definecolor{mycolor1}{rgb}{0.00000,0.44700,0.74100}%
%\definecolor{mycolor2}{rgb}{0.85000,0.32500,0.09800}%
%\definecolor{mycolor3}{rgb}{0.92900,0.69400,0.12500}%
% all in black
\definecolor{mycolor1}{rgb}{0,0,0}%
\definecolor{mycolor2}{rgb}{0,0,0}%
\definecolor{mycolor3}{rgb}{0,0,0}%
\begin{tikzpicture}

\begin{axis}[%
%width=0.75\textwidth,
width=0.6\textwidth,
at={(0.812in,0.573in)},
scale only axis,
xmin=0,
xmax=5,
xlabel style={font=\color{white!15!black}},
xlabel={time $t$},
xtick={0, 1, 2, 3, 4, 5},
ymin=0,
ymax=7,
axis background/.style={fill=white},
axis x line*=bottom,
axis y line*=left,
legend style={at={(0.03,0.97)}, anchor=north west, legend cell align=left, align=left, draw=white!15!black}
]
\addplot [color=mycolor1, dashed, line width=1.2pt]
  table[row sep=crcr]{%
0	1\\
0.125	1.48608530364344\\
0.25	1.90937053770297\\
0.375	2.28681795053251\\
0.5	2.62836232472412\\
0.625	2.94045218747739\\
0.75	3.22757344551866\\
0.875	3.49300597968362\\
1	3.7392392461983\\
1.125	3.96821837816959\\
1.25	4.18149864531897\\
1.375	4.38034700161217\\
1.5	4.56581142372099\\
1.625	4.73876977070432\\
1.75	4.89996513931178\\
1.875	5.05003203301029\\
2	5.18951611228081\\
2.125	5.31888935362826\\
2.25	5.43856185577496\\
2.375	5.54889115175696\\
2.5	5.65018963440833\\
2.625	5.74273053276143\\
2.75	5.82675275948864\\
2.875	5.90246486690547\\
3	5.97004828994225\\
3.125	6.02966001147173\\
3.25	6.08143475364042\\
3.375	6.12548677503282\\
3.5	6.1619113354047\\
3.625	6.19078587570609\\
3.75	6.21217095011474\\
3.875	6.22611093795722\\
4	6.2326345561358\\
4.125	6.23175518652604\\
4.25	6.22347102742673\\
4.375	6.20776507322269\\
4.5	6.18460492170926\\
4.625	6.15394240381301\\
4.75	6.1157130255265\\
4.875	6.06983520724381\\
5	6.01622599875004\\
};
\addlegendentry{$z_s$}

\addplot [color=mycolor2, line width=1.2pt]
  table[row sep=crcr]{%
0	0.5\\
0.125	0.597760945703251\\
0.25	0.65628290457263\\
0.375	0.695754285595634\\
0.5	0.724393509796211\\
0.625	0.746221978578163\\
0.75	0.76345765704472\\
0.875	0.77743183310764\\
1	0.788995651950224\\
1.125	0.798720587969524\\
1.25	0.80700563063967\\
1.375	0.81413837164477\\
1.5	0.820331669243357\\
1.625	0.82574661784661\\
1.75	0.830507466844388\\
1.875	0.834711613542004\\
2	0.838436476947183\\
2.125	0.841744336713501\\
2.25	0.844685809429026\\
2.375	0.847302390986878\\
2.5	0.849628345463429\\
2.625	0.851692128076113\\
2.75	0.8535174702172\\
2.875	0.855124215480672\\
3	0.856528969452336\\
3.125	0.85774560817817\\
3.25	0.858785677896511\\
3.375	0.859658709852482\\
3.5	0.860372467788839\\
3.625	0.860933141140638\\
3.75	0.861345493582463\\
3.875	0.861612974015989\\
4	0.86173779509588\\
4.125	0.861720982794228\\
4.25	0.861562399166032\\
4.375	0.861260739296569\\
4.5	0.860813502299736\\
4.625	0.86021693511579\\
4.75	0.859465946648409\\
4.875	0.858553988410107\\
5	0.857473235401953\\
};
\addlegendentry{$z_f$}

\addplot[const plot, color=mycolor3, dashdotted, line width=1.2pt] table[row sep=crcr] {%
0	4.16500521160618\\
0.125	3.70059991205251\\
0.25	3.35800864710272\\
0.375	3.08763809772788\\
0.5	2.8645290153039\\
0.625	2.67449497484869\\
0.75	2.50875628695749\\
0.875	2.36152637264977\\
1	2.22880180302649\\
1.125	2.10770381649697\\
1.25	1.99609609663653\\
1.375	1.89235106350446\\
1.5	1.7952007110555\\
1.625	1.70363791472119\\
1.75	1.61684909618205\\
1.875	1.53416704744375\\
2	1.45503710277316\\
2.125	1.3789923787439\\
2.25	1.30563531569973\\
2.375	1.23462368581333\\
2.5	1.16565982376809\\
2.625	1.09848221893553\\
2.75	1.03285886193586\\
2.875	0.968581910290269\\
3	0.905463355720588\\
3.125	0.843331458426005\\
3.25	0.782027771792468\\
3.375	0.72140462308847\\
3.5	0.661322945633608\\
3.625	0.601650380082568\\
3.75	0.542259578237133\\
3.875	0.483026654350406\\
4	0.423829736832538\\
4.125	0.364547578706303\\
4.25	0.305058188222017\\
4.375	0.24523744215204\\
4.5	0.184957643704161\\
4.625	0.124085984798275\\
4.75	0.0624828740876031\\
4.875	0.000133757334054843\\
5	nan\\
};
\addlegendentry{$u$}
\end{axis}
\end{tikzpicture}%
		\caption{Numerical solution of  \eqref{eq:origOCP} with multiple shooting.}
		\label{fig:enzymeFull}
	\emp
	\qquad
	\bmp{0.4}
%\definecolor{mycolor1}{rgb}{0.00000,0.44700,0.74100}%
%\definecolor{mycolor2}{rgb}{0.85000,0.32500,0.09800}%
%\definecolor{mycolor3}{rgb}{0.92900,0.69400,0.12500}%
% all in black
\definecolor{mycolor1}{rgb}{0,0,0}%
\definecolor{mycolor2}{rgb}{0,0,0}%
\definecolor{mycolor3}{rgb}{0,0,0}%
\begin{tikzpicture}

\begin{axis}[%
%width=0.75\textwidth,
width=0.6\textwidth,
at={(0.812in,0.573in)},
scale only axis,
xmin=0,
xmax=5,
xlabel style={font=\color{white!15!black}},
xlabel={time $t$},
xtick={0, 1, 2, 3, 4, 5},
ymin=0,
ymax=7,
axis background/.style={fill=white},
axis x line*=bottom,
axis y line*=left,
legend style={at={(0.03,0.97)}, anchor=north west, legend cell align=left, align=left, draw=white!15!black}
]
\addplot [color=mycolor1, dashed, line width=1.2pt]
  table[row sep=crcr]{%
	0	1\\
	0.125	1.4735260711846\\
	0.25	1.88868718683413\\
	0.375	2.26047165895943\\
	0.5	2.59783565202751\\
	0.625	2.90671968648864\\
	0.75	3.19131606202187\\
	0.875	3.45472146088849\\
	1	3.69930428365045\\
	1.125	3.92692603021579\\
	1.25	4.13908210392899\\
	1.375	4.33699534967593\\
	1.5	4.52168045983038\\
	1.625	4.69398967107515\\
	1.75	4.85464602136051\\
	1.875	5.00426808580237\\
	2	5.14338872335056\\
	2.125	5.27246951749778\\
	2.25	5.3919120585253\\
	2.375	5.50206686698652\\
	2.5	5.60324052673283\\
	2.625	5.69570143840293\\
	2.75	5.77968449510125\\
	2.875	5.85539490484225\\
	3	5.92301132891086\\
	3.125	5.98268846484046\\
	3.25	6.03455917274738\\
	3.375	6.07873622123919\\
	3.5	6.11531371193266\\
	3.625	6.14436822830335\\
	3.75	6.16595974409165\\
	3.875	6.18013231804786\\
	4	6.18691459484497\\
	4.125	6.18632012609333\\
	4.25	6.17834752021698\\
	4.375	6.16298042522102\\
	4.5	6.14018734385091\\
	4.625	6.10992127610477\\
	4.75	6.0721191793417\\
	4.875	6.02670123179845\\
	5	5.97358525732899\\
};
\addlegendentry{$z_s$}

\addplot [color=mycolor2, line width=1.2pt]
table[row sep=crcr]{%
	0	0.501\\
	0.125	0.595718835693901\\
	0.25	0.653821983717124\\
	0.375	0.693295907893538\\
	0.5	0.722055119600457\\
	0.625	0.74403077767302\\
	0.75	0.76141145520827\\
	0.875	0.775519073688493\\
	1	0.787202543261746\\
	1.125	0.797033689187292\\
	1.375	0.812628654424306\\
	1.625	0.824376218123491\\
	1.875	0.833451806996329\\
	2.25	0.843552290637942\\
	2.625	0.850650449516082\\
	3.125	0.856788684611201\\
	3.625	0.860029611010479\\
	4.125	0.860846722320506\\
	4.625	0.859351466610351\\
	5	0.856601738833652\\
};
\addlegendentry{$z_f$}

\addplot[const plot, color=mycolor3, dashdotted, line width=1.2pt] table[row sep=crcr] {%
	0	4.15608194376691\\
	0.125	3.70377623669959\\
	0.25	3.36600275790702\\
	0.375	3.09762597158964\\
	0.5	2.87527479530046\\
	0.625	2.68540460441637\\
	0.75	2.51952787564659\\
	0.875	2.37200260170985\\
	1	2.23890128316672\\
	1.125	2.11738674708632\\
	1.25	2.00534609136247\\
	1.375	1.90116510133139\\
	1.5	1.80358338774935\\
	1.625	1.7115980684391\\
	1.75	1.6243977898074\\
	1.875	1.54131634836741\\
	2	1.46179934239018\\
	2.125	1.3853797058718\\
	2.25	1.31165943240622\\
	2.375	1.24029569765314\\
	2.5	1.17099016200001\\
	2.625	1.10348060799549\\
	2.75	1.03753431515255\\
	2.875	0.972942742817627\\
	3	0.909517207690752\\
	3.125	0.84708532369707\\
	3.25	0.785488029443742\\
	3.375	0.72457706975892\\
	3.5	0.664212827655877\\
	3.625	0.604262424761153\\
	3.75	0.544598024011969\\
	3.875	0.485095279784808\\
	4	0.425631888602174\\
	4.125	0.366086198877894\\
	4.25	0.306335841255693\\
	4.375	0.246256342241896\\
	4.5	0.185719683192107\\
	4.625	0.124592764663551\\
	4.75	0.0627357377061921\\
	4.875	0.000124271601594117\\
};
\addlegendentry{$u$}
\end{axis}
\end{tikzpicture}%
		\caption{Numerical solution of  \eqref{eq:liftedOCP} using the ZDP with $m=2$.}
		\label{fig:enzymeZDP}
	\emp
\end{figure}
\begin{figure}[ht]
	\centering
	\bmp{0.4}
%\definecolor{mycolor1}{rgb}{0.00000,0.44700,0.74100}%
%\definecolor{mycolor2}{rgb}{0.85000,0.32500,0.09800}%
%\definecolor{mycolor3}{rgb}{0.92900,0.69400,0.12500}%
% all in black
\definecolor{mycolor1}{rgb}{0,0,0}%
\definecolor{mycolor2}{rgb}{0,0,0}%
\definecolor{mycolor3}{rgb}{0,0,0}%
\begin{tikzpicture}

\begin{axis}[%
%width = 0.7\textwidth,
width=0.6\textwidth,
at={(0.812in,0.573in)},
scale only axis,
xmin=0,
xmax=5,
xlabel style={font=\color{white!15!black}},
xlabel={time $t$},
xtick={0, 1, 2, 3, 4, 5},
ymin=0,
ymax=0.08,
axis background/.style={fill=white},
axis x line*=bottom,
axis y line*=left,
legend style={at={(0.03,0.97)}, anchor=north west, legend cell align=left, align=left, draw=white!15!black}
]
%\addplot [color=mycolor1, line width=1.2pt, mark=x, mark options={solid, mycolor1}]
\addplot [color=mycolor1, line width=1.2pt, dashed]
  table[row sep=crcr]{%
	0	0\\
	0.125	0.0125592324588384\\
	0.25	0.0206833508688469\\
	0.375	0.0263462915730877\\
	0.5	0.0305266726966105\\
	0.625	0.0337325009887515\\
	0.75	0.036257383496797\\
	0.875	0.0382845187951295\\
	1	0.0399349625478482\\
	1.125	0.0412923479538039\\
	1.25	0.042416541389982\\
	1.375	0.0433516519362431\\
	1.5	0.0441309638906082\\
	1.625	0.0447800996291612\\
	1.75	0.0453191179512675\\
	1.875	0.0457639472079201\\
	2	0.0461273889302456\\
	2.125	0.0464198361304842\\
	2.25	0.0466497972496596\\
	2.375	0.0468242847704357\\
	2.5	0.0469491076755002\\
	2.625	0.0470290943584972\\
	2.75	0.0470682643873968\\
	2.875	0.0470699620632216\\
	3	0.0470369610313899\\
	3.125	0.0469715466312683\\
	3.25	0.0468755808930421\\
	3.375	0.0467505537936352\\
	3.5	0.0465976234720413\\
	3.625	0.046417647402734\\
	3.75	0.0462112060230897\\
	3.875	0.0459786199093584\\
	4	0.0457199612908274\\
	4.125	0.0454350604327152\\
	4.25	0.0451235072097473\\
	4.375	0.0447846480016638\\
	4.5	0.0444175778583453\\
	4.625	0.0440211277082376\\
	4.75	0.0435938461847973\\
	4.875	0.0431339754453619\\
	5	0.0426407414210521\\
};
\addlegendentry{$z_s$}

%\addplot [color=mycolor2, line width=1.2pt, mark=o, mark options={solid, mycolor2}]
\addplot [color=mycolor2, line width=1.2pt]
  table[row sep=crcr]{%
	0	0.00100000000000033\\
	0.125	0.00204211000934862\\
	0.25	0.00246092085550664\\
	0.375	0.00245837770209612\\
	0.5	0.00233839019575477\\
	0.625	0.00219120090514302\\
	0.75	0.00204620183644977\\
	0.875	0.00191275941914704\\
	1	0.00179310868847882\\
	1.125	0.00168689878223294\\
	1.25	0.001592909024553\\
	1.375	0.00150971722046389\\
	1.5	0.00143595575085964\\
	1.625	0.00137039972312003\\
	1.75	0.00131198625538698\\
	1.875	0.0012598065456757\\
	2	0.00121308852624047\\
	2.125	0.00117117764280117\\
	2.25	0.00113351879108414\\
	2.375	0.00109964045494859\\
	2.5	0.00106914122901713\\
	2.625	0.00104167856003112\\
	2.75	0.0010169594224223\\
	2.875	0.000994732625302142\\
	3	0.00097478247938465\\
	3.125	0.000956923566969081\\
	3.25	0.000940996437877928\\
	3.375	0.000926864046387088\\
	3.5	0.000914408803042832\\
	3.625	0.000903530130160313\\
	3.75	0.00089414243225594\\
	3.875	0.000886173409044844\\
	4	0.000879562651777022\\
	4.125	0.000874260473722188\\
	4.25	0.000870226933302298\\
	4.375	0.000867431013745978\\
	4.5	0.000865849926327655\\
	4.625	0.000865468505439004\\
	4.75	0.000866278662198461\\
	4.875	0.000868278859034177\\
	5	0.000871496568300323\\
};
\addlegendentry{$z_f$}

%\addplot [color=mycolor3, line width=1.2pt, mark=square, mark options={solid, mycolor3}]
\addplot [color=mycolor3, line width=1.2pt, dashdotted]
  table[row sep=crcr]{%
	0	0.00892326783926478\\
	0.125	0.00317632464708417\\
	0.25	0.0079941108042938\\
	0.375	0.00998787386175959\\
	0.5	0.0107457799965536\\
	0.625	0.0109096295676849\\
	0.75	0.0107715886890993\\
	0.875	0.0104762290600817\\
	1	0.0100994801402354\\
	1.125	0.00968293058934844\\
	1.25	0.00924999472594124\\
	1.375	0.00881403782692658\\
	1.5	0.0083826766938504\\
	1.625	0.00796015371791281\\
	1.75	0.00754869362535082\\
	1.875	0.00714930092365851\\
	2	0.00676223961702505\\
	2.125	0.00638732712789025\\
	2.25	0.00602411670649428\\
	2.375	0.00567201183980881\\
	2.5	0.00533033823191786\\
	2.625	0.00499838905995542\\
	2.75	0.0046754532166835\\
	2.875	0.0043608325273583\\
	3	0.00405385197016361\\
	3.125	0.00375386527106514\\
	3.25	0.00346025765127411\\
	3.375	0.00317244667045014\\
	3.5	0.00288988202226914\\
	3.625	0.00261204467858533\\
	3.75	0.00233844577483566\\
	3.875	0.00206862543440156\\
	4	0.00180215176963561\\
	4.125	0.00153862017159057\\
	4.25	0.001277653033676\\
	4.375	0.00101890008985617\\
	4.5	0.000762039487946176\\
	4.625	0.000506779865275142\\
	4.75	0.000252863618588783\\
	4.875	9.4857324608455e-06\\
};
\addlegendentry{$u$}
\end{axis}
\end{tikzpicture}%
		\caption{absolute error between full \eqref{eq:origOCP} and lifted solution \eqref{eq:liftedOCP}.}
		\label{fig:enzymeDiffAbs}
	\emp
	\qquad
	\bmp{0.4}
%\definecolor{mycolor1}{rgb}{0.00000,0.44700,0.74100}%
%\definecolor{mycolor2}{rgb}{0.85000,0.32500,0.09800}%
%\definecolor{mycolor3}{rgb}{0.92900,0.69400,0.12500}%
% all in black
\definecolor{mycolor1}{rgb}{0,0,0}%
\definecolor{mycolor2}{rgb}{0,0,0}%
\definecolor{mycolor3}{rgb}{0,0,0}%
\begin{tikzpicture}

\begin{axis}[%
%width = 0.7\textwidth,
width=0.6\textwidth,
at={(0.812in,0.573in)},
scale only axis,
xmin=0,
xmax=5,
xlabel style={font=\color{white!15!black}},
xlabel={time $t$},
xtick={0, 1, 2, 3, 4, 5},
ymode=log,
ymin=0.000858132851115336,
ymax=1,
yminorticks=true,
axis background/.style={fill=white},
legend style={at={(0.03,0.97)}, anchor=north west, legend cell align=left, align=left, draw=white!15!black}
]
%\addplot [color=mycolor1, line width=1.2pt, mark=x, mark options={solid, mycolor1}]
\addplot [color=mycolor1, line width=1.2pt, dashed]
table[row sep=crcr]{%
0	0\\
0.125	0.00845102412982989\\
0.25	0.0108324758942293\\
0.375	0.0115209152136607\\
0.5	0.0116143312422228\\
0.625	0.0114718872134779\\
0.75	0.0112336545378364\\
0.875	0.0109603593021433\\
1	0.0106799959886593\\
1.125	0.0104057928207963\\
1.25	0.0101438900037225\\
1.375	0.00989688300881144\\
1.5	0.00966555502132502\\
1.625	0.00944976150123578\\
1.75	0.00924889613908174\\
1.875	0.00906214124999364\\
2	0.0088886034129116\\
2.125	0.00872738644526434\\
2.25	0.00857762953192502\\
2.375	0.00843852577056024\\
2.5	0.00830932968832452\\
2.625	0.00818935860374059\\
2.75	0.00807799064058735\\
2.875	0.00797466101765118\\
3	0.00787885755056092\\
3.125	0.00779011589742993\\
3.25	0.0077080148395584\\
3.375	0.00763217174491707\\
3.5	0.0075622382755379\\
3.625	0.0074978963474144\\
3.75	0.00743885431887938\\
3.875	0.00738484336209525\\
4	0.00733561395677355\\
4.125	0.00729093243186399\\
4.25	0.00725057746701928\\
4.375	0.00721433644888908\\
4.5	0.00718200155549725\\
4.625	0.00715336541251744\\
4.75	0.00712821612500826\\
4.875	0.00710633143319487\\
5	0.00708898228073414\\
};
\addlegendentry{$z_s$}

%\addplot [color=mycolor2, line width=1.2pt, mark=o, mark options={solid, mycolor2}]
\addplot [color=mycolor2, line width=1.2pt]
  table[row sep=crcr]{%
0	0.002\\
0.125	0.00341656263690394\\
0.25	0.00374995306795107\\
0.375	0.00353350581786531\\
0.5	0.00322814049748435\\
0.625	0.00293644793438775\\
0.75	0.00268021988094404\\
0.875	0.00246039061531837\\
1	0.00227267524879446\\
1.125	0.00211202481288943\\
1.25	0.00197387149563055\\
1.375	0.00185439191450433\\
1.5	0.00175047313222488\\
1.625	0.00165960258218528\\
1.75	0.00157975302217722\\
1.875	0.00150928303696465\\
2	0.00144685655166974\\
2.125	0.00139137947221845\\
2.25	0.00134195032573983\\
2.375	0.00129782189632049\\
2.5	0.00125837135105334\\
2.625	0.0012230768829029\\
2.75	0.00119149935453418\\
2.875	0.00116326779023823\\
3	0.00113806785104686\\
3.125	0.00111563261297638\\
3.25	0.00109573517405053\\
3.375	0.00107818268423658\\
3.5	0.00106281151438334\\
3.625	0.00104948333575608\\
3.75	0.00103808193458004\\
3.875	0.00102851062473199\\
4	0.00102069015129434\\
4.125	0.00101455700027815\\
4.25	0.00101006204701625\\
4.375	0.0010071694886755\\
4.5	0.00100585601569976\\
4.625	0.00100611018321919\\
4.75	0.00100793194655363\\
4.875	0.00101133232479382\\
5	0.00101654889816607\\
};
\addlegendentry{$z_f$}

%\addplot [color=mycolor3, line width=1.2pt, mark=square, mark options={solid, mycolor3}]
\addplot [color=mycolor3, line width=1.2pt, dashdotted]
  table[row sep=crcr]{%
0	0.00214186459930352\\
0.125	0.000858132851115336\\
0.25	0.00238050271531981\\
0.375	0.00323472798142185\\
0.5	0.0037512817461019\\
0.625	0.00407910629429131\\
0.75	0.00429357547363723\\
0.875	0.00443619468304661\\
1	0.00453133685236264\\
1.125	0.00459405553921803\\
1.25	0.00463403356566157\\
1.375	0.00465770953681928\\
1.5	0.00466948430297302\\
1.625	0.00467243585704613\\
1.75	0.00466875885007078\\
1.875	0.00466004347531611\\
2	0.00464745733394157\\
2.125	0.0046318668959002\\
2.25	0.00461392033279951\\
2.375	0.00459410505453301\\
2.5	0.00457278832342777\\
2.625	0.00455024632697008\\
2.75	0.00452668523838189\\
2.875	0.00450225662289816\\
3	0.00447706878772529\\
3.125	0.0044511951827636\\
3.25	0.00442468061462351\\
3.375	0.00439754581718356\\
3.5	0.00436979076038354\\
3.625	0.00434139696798575\\
3.75	0.00431232903459069\\
3.875	0.00428253547261774\\
4	0.00425194897441378\\
4.125	0.00422048613823067\\
4.25	0.00418804667599807\\
4.375	0.0041545120925008\\
4.5	0.00411974379416276\\
4.625	0.00408358050392728\\
4.75	0.00404583357698278\\
4.875	0.547278354530861\\
};
\addlegendentry{$u$}
\end{axis}
\end{tikzpicture}%
		\caption{relative error between full \eqref{eq:origOCP} and lifted solution \eqref{eq:liftedOCP}.}
		\label{fig:enzymeDiffRel}
	\emp
\end{figure}

However, there is a big difference in the runtimes which are listed in Table \ref{tab:runtimeEnzyme}. Whereas the NLP resulting from the original system can be solved with CasADi/Ipopt in 0.1688 seconds, our proposed approach using ZDP of order $m=2$ only needs 0.0224 seconds. This gives a significant speed up of factor 7.5. The order reduced OCP \eqref{eq:reducedOCP} takes 0.5272 seconds. It has the smallest number of optimization variables, but the nested root finding problem for $h_\epsilon(\cdot)$ in each interval causes this longer runtime. 
As already mentioned in Rehberg\cite{Rehberg2013}, the order reduced OCP \eqref{eq:reducedOCP} does not need to improve the runtime, although less optimization variables are present. 

%%% PLACEHOLDER FOR TABLE 1

For other values of  $\epsilon \in (0,0.1]$ and the choice of $z_s(0)$ almost the same runtimes are obtained for both the full problem and 
the lifted problem.

\subsection{CSTR Example} 
\begin{figure}[htp]		
	\centering
	\def\svgwidth{0.3\textwidth}
	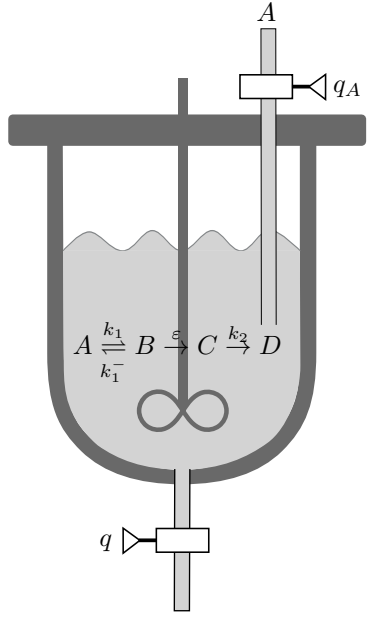
	\caption{Visualization of the CSTR}
	\label{fig:cstr}
\end{figure}

Inspired by an example of Cvejn \cite{Cvejn2007}, we consider a small continuous stirred-tank reactor (CSTR) in which the chemical reactions 
$$A \overset{k_1}{\underset{k_1^-}{\rightleftharpoons}} B\overset{\epsilon}{\rightarrow}C\overset{k_2}{\rightarrow} D$$
take place for given reaction rates $k_1,k_1^-, k_2$ and $\epsilon$. For simplicity, we assume that the tank is perfectly mixed and the reaction is isothermic. Figure \ref{fig:cstr} illustrates the CSTR. 
Input for the reaction is only species $A$ with a given (constant) concentration $c_{A_{\text{in}}}$. This input can be controlled via a feed 
$q_A \in [0,q_{A_\text{max}}]$. Output of the reactor is a mixture of species $A,B,C$ and $D$ and can also be controlled with the outlet
$q \in [0,q_{\text{max}}]$. A possible scenario is that the number of moles of species $B$ should be maximized in the output of the CSTR.
Since the volume of the tank reactor is bounded, this leads to an optimal control problem with constraints (see \eqref{eq:cstrOCP}).
For brevity, we omit the arguments of the functions in the ODE part of the OCP (i.e. the dependency on $t$).

\begin{subequations}\label{eq:cstrOCP}
	\begin{align}
		&\underset{c(\cdot),V(\cdot),q(\cdot),q_A(\cdot)}{\text{minimize}}	 &  &\int_{0}^{500} -0.1 \, q \cdot c_B + q^2 +q_A^2 \dif t & &\\
		& \text{such that} &\dot{c_A} 		&=-k_1 c_A	+ k_1^- c_B + \frac{q_A}{V}\left(c_{A_{\text{in}}}-c_A\right),   & & \\
		&						 &\dot{c_B} 		&=  k_1 c_A	- (k_1^- +\epsilon)c_B	- \frac{q_A}{V}c_B ,& &\\
		&						 &\dot{c_C} 		&=  \epsilon c_B - k_2 c_C - \frac{q_A}{V}c_C,& &\\
		&						 &\dot{c_D} 		&=  k_2 c_C	- \frac{q_A}{V}c_D,& &\\
		&						 &\dot{V} 			&=  q_A-q,& &\\
		&						 & c^*	&= c(0),& & \\
		& 						 & V^* &=V(0) =V(500),
	\end{align}
\end{subequations}

where $c(t) = \left(c_A(t),c_B(t),c_C(t),c_D(t)\right)^T$ and $c^* = \left(c_A^*,c_B^*,c_C^*,c_D^*\right)^T$. 	
The values of the constants in OCP \eqref{eq:cstrOCP} can be taken from Table \ref{table:cstrConstants}.

%%% PLACEHOLDER FOR TABLE 2

Note that the ODE of this OCP is not in singularly perturbed form. However, the idea of lifting in combination with the ZDP method is still
applicable, because its dynamics are given by a multi-scale ODE.
The choice of reaction progress variables is not obvious and has combinatorial complexity. We choose $B,D$ and the volume $V$ mainly motivated by 
the reaction rates.

For $N=140$ intervals IPOPT converges to a point of local infeasibility for the nonlinear program resulting from the lifted OCP using multiple shooting 
and ZDP of order $2$. Without the constraint $c_A(0) = c_A^*$, it converges to the solution depicted in 
Figures \ref{fig:cstrFullStates},\ref{fig:cstrZDPStates},\ref{fig:cstrFullControls},\ref{fig:cstrZDPControls},\ref{fig:cstrFullVolume}, and \ref{fig:cstrZDPVolume}.
For the plots of the lifted OCP, we took the optimal controls of the numerical solution and integrated the other states forward in time with an implicit Radau II-A integrator.

\begin{figure}[ht]
	\centering
	\bmp{0.4}
		\input{cstrFullN140States}
		\caption{``full'' species $c_A,c_B,c_C$ and $c_D$ of  \eqref{eq:origOCP} solved by multiple shooting.}
		\label{fig:cstrFullStates}
	\emp
	\qquad
	\bmp{0.4}
		\input{cstrZDPN140States}
		\caption{``lifted'' species $c_A,c_B,c_C$ and $c_D$ of  \eqref{eq:liftedOCP} using the ZDP with $m=2$. Plotted values are obtained by forward integration with the calculated ``lifted'' controls.}
		\label{fig:cstrZDPStates}
	\emp
\end{figure}
\begin{figure}[ht]
 \centering
	\bmp{0.4}
%\definecolor{mycolor1}{rgb}{0.00000,0.44700,0.74100}%
%\definecolor{mycolor2}{rgb}{0.85000,0.32500,0.09800}%
% all in black
\definecolor{mycolor1}{rgb}{0,0,0}%
\definecolor{mycolor2}{rgb}{0,0,0}%
\begin{tikzpicture}

\begin{axis}[%
%width = 0.75\textwidth,
width=0.6\textwidth,
at={(0.812in,0.573in)},
scale only axis,
xmin=0,
xmax=500,
xlabel style={font=\color{white!15!black}},
xlabel={time $t$},
xtick={0,100,200,300,400,500},
ymin=0,
ymax=0.0015,
axis background/.style={fill=white},
axis x line*=bottom,
axis y line*=left,
legend style={at={(0.97,0.5)}, anchor=east, legend cell align=left, align=left, draw=white!15!black}
]
\addplot[color=mycolor1, line width=1.2pt] table[row sep=crcr] {%
0	0.001\\
3.57142857142857	0.001\\
7.14285714285714	0.001\\
10.7142857142857	0.001\\
14.2857142857143	0.001\\
17.8571428571429	0.001\\
21.4285714285714	0.001\\
25	0.001\\
28.5714285714286	0.001\\
32.1428571428571	0.001\\
35.7142857142857	0.001\\
39.2857142857143	0.001\\
42.8571428571429	0.001\\
46.4285714285714	0.001\\
50	0.001\\
53.5714285714286	0.001\\
57.1428571428571	0.001\\
60.7142857142857	0.001\\
64.2857142857143	0.001\\
67.8571428571429	0.001\\
71.4285714285714	0.001\\
75	0.001\\
78.5714285714286	0.001\\
82.1428571428571	0.001\\
85.7142857142857	0.001\\
89.2857142857143	0.001\\
92.8571428571429	0.001\\
96.4285714285714	0.001\\
100	0.001\\
103.571428571429	0.001\\
107.142857142857	0.001\\
110.714285714286	0.001\\
114.285714285714	0.001\\
117.857142857143	0.001\\
121.428571428571	0.001\\
125	0.001\\
128.571428571429	0.001\\
132.142857142857	0.001\\
135.714285714286	0.001\\
139.285714285714	0.001\\
142.857142857143	0.001\\
146.428571428571	0.001\\
150	0.001\\
153.571428571429	0.001\\
157.142857142857	0.001\\
160.714285714286	0.001\\
164.285714285714	0.001\\
167.857142857143	0.001\\
171.428571428571	0.001\\
175	0.001\\
178.571428571429	0.001\\
182.142857142857	0.001\\
185.714285714286	0.001\\
189.285714285714	0.001\\
192.857142857143	0.001\\
196.428571428571	0.001\\
200	0.001\\
203.571428571429	0.001\\
207.142857142857	0.001\\
210.714285714286	0.001\\
214.285714285714	0.001\\
217.857142857143	0.001\\
221.428571428571	0.001\\
225	0.001\\
228.571428571429	0.001\\
232.142857142857	0.001\\
235.714285714286	0.001\\
239.285714285714	0.001\\
242.857142857143	0.001\\
246.428571428571	0.001\\
250	0.001\\
253.571428571429	0.001\\
257.142857142857	0.001\\
260.714285714286	0.001\\
264.285714285714	0.001\\
267.857142857143	0.001\\
271.428571428571	0.001\\
275	0.001\\
278.571428571429	0.001\\
282.142857142857	0.001\\
285.714285714286	0.001\\
289.285714285714	0.001\\
292.857142857143	0.001\\
296.428571428571	0.001\\
300	0.001\\
303.571428571429	0.001\\
307.142857142857	0.001\\
310.714285714286	0.001\\
314.285714285714	0.001\\
317.857142857143	0.001\\
321.428571428571	0.001\\
325	0.001\\
328.571428571429	0.001\\
332.142857142857	0.001\\
335.714285714286	0.001\\
339.285714285714	0.001\\
342.857142857143	0.001\\
346.428571428571	0.001\\
350	0.001\\
353.571428571429	0.001\\
357.142857142857	0.001\\
360.714285714286	0.001\\
364.285714285714	0.001\\
367.857142857143	0.001\\
371.428571428571	0.001\\
375	0.001\\
378.571428571429	0.001\\
382.142857142857	0.001\\
385.714285714286	0.001\\
389.285714285714	0.001\\
392.857142857143	0.001\\
396.428571428571	0.001\\
400	0.001\\
403.571428571429	0.001\\
407.142857142857	0.001\\
410.714285714286	0.001\\
414.285714285714	0.001\\
417.857142857143	0.001\\
421.428571428571	0.001\\
425	0.001\\
428.571428571429	0.001\\
432.142857142857	0.001\\
435.714285714286	0.001\\
439.285714285714	0.001\\
442.857142857143	0.001\\
446.428571428571	0.001\\
450	0.001\\
453.571428571429	0.001\\
457.142857142857	0.001\\
460.714285714286	0.001\\
464.285714285714	0.001\\
467.857142857143	0.001\\
471.428571428571	0.001\\
475	0.001\\
478.571428571429	0.001\\
482.142857142857	0.001\\
485.714285714286	0.001\\
489.285714285714	0.001\\
492.857142857143	0.001\\
496.428571428571	0.001\\
500	nan\\
};
\addlegendentry{$q_{\textrm{A}}$}

\addplot[color=mycolor2, dashed, line width=1.2pt] table[row sep=crcr] {%
0	8.28361571172636e-10\\
3.57142857142857	1.02110765285126e-08\\
7.14285714285714	2.3711232885265e-08\\
10.7142857142857	4.24883634142645e-08\\
14.2857142857143	6.83301262399999e-08\\
17.8571428571429	1.04026795258484e-07\\
21.4285714285714	1.54174472980856e-07\\
25	2.26901591787526e-07\\
28.5714285714286	3.38044550723608e-07\\
32.1428571428571	5.23070709262872e-07\\
35.7142857142857	8.80933594720702e-07\\
39.2857142857143	1.82784964814758e-06\\
42.8571428571429	9.24773095427989e-06\\
46.4285714285714	0.000103335779727676\\
50	0.000206886558607564\\
53.5714285714286	0.000293540589016337\\
57.1428571428571	0.000367025714377562\\
60.7142857142857	0.000430314067000178\\
64.2857142857143	0.00048554451750344\\
67.8571428571429	0.000534279222704101\\
71.4285714285714	0.000577685772350771\\
75	0.000616656122002178\\
78.5714285714286	0.000651884984366835\\
82.1428571428571	0.000683922706559825\\
85.7142857142857	0.000713211812662107\\
89.2857142857143	0.000740112836119553\\
92.8571428571429	0.000764922959591311\\
96.4285714285714	0.000787889715191024\\
100	0.000809221222283035\\
103.571428571429	0.000829093952865895\\
107.142857142857	0.000847658701414565\\
110.714285714286	0.000865045230761365\\
114.285714285714	0.000881365927985124\\
117.857142857143	0.000896718710554867\\
121.428571428571	0.000911189358004186\\
125	0.000924853398723485\\
128.571428571429	0.00093777764881017\\
132.142857142857	0.000950021476307856\\
135.714285714286	0.000961637846921451\\
139.285714285714	0.000972674194519371\\
142.857142857143	0.000983173150082007\\
146.428571428571	0.000993173155650569\\
150	0.00100270898415311\\
153.571428571429	0.00101181218190317\\
157.142857142857	0.00102051144703539\\
160.714285714286	0.00102883295477475\\
164.285714285714	0.00103680063825039\\
167.857142857143	0.00104443643197309\\
171.428571428571	0.00105176048382044\\
175	0.00105879134034379\\
178.571428571429	0.00106554610932851\\
182.142857142857	0.00107204060290559\\
185.714285714286	0.00107828946395948\\
189.285714285714	0.00108430627811751\\
192.857142857143	0.0010901036732375\\
196.428571428571	0.00109569340807182\\
200	0.00110108645140063\\
203.571428571429	0.00110629305294466\\
207.142857142857	0.00111132280695753\\
210.714285714286	0.00111618470942109\\
214.285714285714	0.00112088720960339\\
217.857142857143	0.00112543825663407\\
221.428571428571	0.00112984534162548\\
225	0.00113411553590959\\
228.571428571429	0.00113825552576771\\
232.142857142857	0.00114227164407032\\
235.714285714286	0.00114616989916134\\
239.285714285714	0.00114995600127488\\
242.857142857143	0.00115363538676225\\
246.428571428571	0.00115721324037367\\
250	0.00116069451577055\\
253.571428571429	0.0011640839545138\\
257.142857142857	0.00116738610363641\\
260.714285714286	0.00117060533199461\\
264.285714285714	0.00117374584553548\\
267.857142857143	0.00117681170156342\\
271.428571428571	0.00117980682217031\\
275	0.0011827350069121\\
278.571428571429	0.00118559994481962\\
282.142857142857	0.00118840522585862\\
285.714285714286	0.0011911543518823\\
289.285714285714	0.00119385074722011\\
292.857142857143	0.00119649776889596\\
296.428571428571	0.00119909871663842\\
300	0.00120165684268422\\
303.571428571429	0.00120417536149345\\
307.142857142857	0.00120665745940744\\
310.714285714286	0.00120910630436946\\
314.285714285714	0.0012115250557503\\
317.857142857143	0.0012139168743664\\
321.428571428571	0.00121628493278984\\
325	0.00121863242603396\\
328.571428571429	0.00122096258268576\\
332.142857142857	0.00122327867665497\\
335.714285714286	0.00122558403959367\\
339.285714285714	0.00122788207417744\\
342.857142857143	0.00123017626837879\\
346.428571428571	0.00123247021088517\\
350	0.00123476760790553\\
353.571428571429	0.00123707230156435\\
357.142857142857	0.00123938829017407\\
360.714285714286	0.00124171975065351\\
364.285714285714	0.00124407106349686\\
367.857142857143	0.00124644684072492\\
371.428571428571	0.00124885195729297\\
375	0.00125129158656301\\
378.571428571429	0.00125377124059937\\
382.142857142857	0.00125629681609007\\
385.714285714286	0.00125887464700921\\
389.285714285714	0.00126151156521965\\
392.857142857143	0.00126421497060116\\
396.428571428571	0.00126699291259479\\
400	0.0012698541855394\\
403.571428571429	0.00127280844063617\\
407.142857142857	0.00127586631829134\\
410.714285714286	0.00127903960535608\\
414.285714285714	0.00128234142317783\\
417.857142857143	0.00128578645389542\\
421.428571428571	0.00128939121465023\\
425	0.0012931743923341\\
428.571428571429	0.00129715725535952\\
432.142857142857	0.00130136416454053\\
435.714285714286	0.00130582321259767\\
439.285714285714	0.00131056703249986\\
442.857142857143	0.00131563382998105\\
446.428571428571	0.00132106871747169\\
450	0.00132692545856245\\
453.571428571429	0.00133326877971676\\
457.142857142857	0.00134017747722017\\
460.714285714286	0.00134774865534344\\
464.285714285714	0.00135610359561471\\
467.857142857143	0.00136539599972438\\
471.428571428571	0.00137582367868309\\
475	0.0013876450679018\\
478.571428571429	0.00140120149308849\\
482.142857142857	0.00141694085625963\\
485.714285714286	0.0014354046904622\\
489.285714285714	0.00145691118236246\\
492.857142857143	0.00147906530915676\\
496.428571428571	0.00149209188676547\\
500	nan\\
};
\addlegendentry{$q$}
\end{axis}
\end{tikzpicture}%
		\caption{``full'' controls $q_A$ and $q$ of  \eqref{eq:origOCP} solved by multiple shooting.}
		\label{fig:cstrFullControls}
	\emp
	\qquad
	\bmp{0.4}
%\definecolor{mycolor1}{rgb}{0.00000,0.44700,0.74100}%
%\definecolor{mycolor2}{rgb}{0.85000,0.32500,0.09800}%
% all in black
\definecolor{mycolor1}{rgb}{0,0,0}%
\definecolor{mycolor2}{rgb}{0,0,0}%
\begin{tikzpicture}

\begin{axis}[%
%width = 0.75\textwidth,
width=0.6\textwidth,
at={(0.812in,0.573in)},
scale only axis,
unbounded coords=jump,
xmin=0,
xmax=500,
xlabel style={font=\color{white!15!black}},
xlabel={time $t$},
xtick={0,100,200,300,400,500},
ymin=0,
ymax=0.0015,
axis background/.style={fill=white},
axis x line*=bottom,
axis y line*=left,
legend style={at={(0.97,0.5)}, anchor=east, legend cell align=left, align=left, draw=white!15!black}
]
\addplot [color=mycolor1, line width=1.2pt]
table[row sep=crcr]{%
	0	0.001\\
	3.57142857142857	0.001\\
	7.14285714285714	0.001\\
	10.7142857142857	0.001\\
	14.2857142857143	0.001\\
	17.8571428571429	0.001\\
	21.4285714285714	0.001\\
	25	0.001\\
	28.5714285714286	0.001\\
	32.1428571428571	0.001\\
	35.7142857142857	0.001\\
	39.2857142857143	0.001\\
	42.8571428571429	0.001\\
	46.4285714285714	0.001\\
	50	0.001\\
	53.5714285714286	0.001\\
	57.1428571428571	0.001\\
	60.7142857142857	0.001\\
	64.2857142857143	0.001\\
	67.8571428571429	0.001\\
	71.4285714285714	0.001\\
	75	0.001\\
	78.5714285714286	0.001\\
	82.1428571428571	0.001\\
	85.7142857142857	0.001\\
	89.2857142857143	0.001\\
	92.8571428571429	0.001\\
	96.4285714285714	0.001\\
	100	0.001\\
	103.571428571429	0.001\\
	107.142857142857	0.001\\
	110.714285714286	0.001\\
	114.285714285714	0.001\\
	117.857142857143	0.001\\
	121.428571428571	0.001\\
	125	0.001\\
	128.571428571429	0.001\\
	132.142857142857	0.001\\
	135.714285714286	0.001\\
	139.285714285714	0.001\\
	142.857142857143	0.001\\
	146.428571428571	0.001\\
	150	0.001\\
	153.571428571429	0.001\\
	157.142857142857	0.001\\
	160.714285714286	0.001\\
	164.285714285714	0.001\\
	167.857142857143	0.001\\
	171.428571428571	0.001\\
	175	0.001\\
	178.571428571429	0.001\\
	182.142857142857	0.001\\
	185.714285714286	0.001\\
	189.285714285714	0.001\\
	192.857142857143	0.001\\
	196.428571428571	0.001\\
	200	0.001\\
	203.571428571429	0.001\\
	207.142857142857	0.001\\
	210.714285714286	0.001\\
	214.285714285714	0.001\\
	217.857142857143	0.001\\
	221.428571428571	0.001\\
	225	0.001\\
	228.571428571429	0.001\\
	232.142857142857	0.001\\
	235.714285714286	0.001\\
	239.285714285714	0.001\\
	242.857142857143	0.001\\
	246.428571428571	0.001\\
	250	0.001\\
	253.571428571429	0.001\\
	257.142857142857	0.001\\
	260.714285714286	0.001\\
	264.285714285714	0.001\\
	267.857142857143	0.001\\
	271.428571428571	0.001\\
	275	0.001\\
	278.571428571429	0.001\\
	282.142857142857	0.001\\
	285.714285714286	0.001\\
	289.285714285714	0.001\\
	292.857142857143	0.001\\
	296.428571428571	0.001\\
	300	0.001\\
	303.571428571429	0.001\\
	307.142857142857	0.001\\
	310.714285714286	0.001\\
	314.285714285714	0.001\\
	317.857142857143	0.001\\
	321.428571428571	0.001\\
	325	0.001\\
	328.571428571429	0.001\\
	332.142857142857	0.001\\
	335.714285714286	0.001\\
	339.285714285714	0.001\\
	342.857142857143	0.001\\
	346.428571428571	0.001\\
	350	0.001\\
	353.571428571429	0.001\\
	357.142857142857	0.001\\
	360.714285714286	0.001\\
	364.285714285714	0.001\\
	367.857142857143	0.001\\
	371.428571428571	0.001\\
	375	0.001\\
	378.571428571429	0.001\\
	382.142857142857	0.001\\
	385.714285714286	0.001\\
	389.285714285714	0.001\\
	392.857142857143	0.001\\
	396.428571428571	0.001\\
	400	0.001\\
	403.571428571429	0.001\\
	407.142857142857	0.001\\
	410.714285714286	0.001\\
	414.285714285714	0.001\\
	417.857142857143	0.001\\
	421.428571428571	0.001\\
	425	0.001\\
	428.571428571429	0.001\\
	432.142857142857	0.001\\
	435.714285714286	0.001\\
	439.285714285714	0.001\\
	442.857142857143	0.001\\
	446.428571428571	0.001\\
	450	0.001\\
	453.571428571429	0.001\\
	457.142857142857	0.001\\
	460.714285714286	0.001\\
	464.285714285714	0.001\\
	467.857142857143	0.001\\
	471.428571428571	0.001\\
	475	0.001\\
	478.571428571429	0.001\\
	482.142857142857	0.001\\
	485.714285714286	0.001\\
	489.285714285714	0.001\\
	492.857142857143	0.001\\
	496.428571428571	0.001\\
	500	nan\\
};
\addlegendentry{$q_{\textrm{A}}$}

\addplot [color=mycolor2, dashed, line width=1.2pt]
table[row sep=crcr]{%
	0	0\\
	3.57142857142857	0.000599365588076857\\
	7.14285714285714	0.000633247129711553\\
	10.7142857142857	0.000660322304226276\\
	14.2857142857143	0.000682754529829665\\
	17.8571428571429	0.00070183831206047\\
	21.4285714285714	0.000718404673070056\\
	25	0.000733016396115747\\
	28.5714285714286	0.000746070960209104\\
	32.1428571428571	0.000757858766846888\\
	35.7142857142857	0.00076859798781413\\
	39.2857142857143	0.000778456407140627\\
	42.8571428571429	0.000787565648122101\\
	46.4285714285714	0.00079603074803544\\
	50	0.000803936785441267\\
	53.5714285714286	0.000811353578544122\\
	57.1428571428571	0.000818339089694717\\
	60.7142857142857	0.000824941937340353\\
	64.2857142857143	0.000831203282369011\\
	67.8571428571429	0.000837158264328331\\
	71.4285714285714	0.000842837111407608\\
	75	0.000848266007448254\\
	78.5714285714286	0.000853467777781153\\
	82.1428571428571	0.000858462436794536\\
	85.7142857142857	0.000863267629135533\\
	89.2857142857143	0.00086789898883828\\
	92.8571428571429	0.000872370432948659\\
	96.4285714285714	0.000876694404635961\\
	100	0.000880882074848205\\
	103.571428571429	0.000884943511113029\\
	107.142857142857	0.000888887819599498\\
	110.714285714286	0.000892723265257157\\
	114.285714285714	0.000896457374414254\\
	117.857142857143	0.000900097021756918\\
	121.428571428571	0.000903648505911696\\
	125	0.000907117613705768\\
	128.571428571429	0.000910509677032602\\
	132.142857142857	0.000913829620989815\\
	135.714285714286	0.00091708200701784\\
	139.285714285714	0.000920271070128738\\
	142.857142857143	0.000923400751703444\\
	146.428571428571	0.000926474728603488\\
	150	0.000929496439110994\\
	153.571428571429	0.000932469105185562\\
	157.142857142857	0.000935395753274963\\
	160.714285714286	0.00093827923234485\\
	164.285714285714	0.000941122230149104\\
	167.857142857143	0.000943927287801462\\
	171.428571428571	0.000946696813224307\\
	175	0.000949433092749205\\
	178.571428571429	0.000952138302183467\\
	182.142857142857	0.000954814516582112\\
	185.714285714286	0.000957463719132344\\
	189.285714285714	0.000960087809721653\\
	192.857142857143	0.00096268861186118\\
	196.428571428571	0.00096526788015419\\
	200	0.000967827306533193\\
	203.571428571429	0.000970368526302009\\
	207.142857142857	0.000972893123400836\\
	210.714285714286	0.000975402635832871\\
	214.285714285714	0.000977898560428614\\
	217.857142857143	0.000980382357409313\\
	221.428571428571	0.000982855454706308\\
	225	0.000985319252065758\\
	228.571428571429	0.000987775124875992\\
	232.142857142857	0.000990224428203349\\
	235.714285714286	0.000992668500257536\\
	239.285714285714	0.000995108665537045\\
	242.857142857143	0.000997546239099431\\
	246.428571428571	0.000999982529216687\\
	250	0.00100241884127179\\
	253.571428571429	0.00100485648034299\\
	257.142857142857	0.00100729675553349\\
	260.714285714286	0.00100974098223155\\
	264.285714285714	0.00101219048658115\\
	267.857142857143	0.00101464660776724\\
	271.428571428571	0.00101711070300302\\
	275	0.00101958414955723\\
	278.571428571429	0.0010220683497923\\
	282.142857142857	0.00102456473467098\\
	285.714285714286	0.00102707476771172\\
	289.285714285714	0.00102959994969021\\
	292.857142857143	0.00103214182299253\\
	296.428571428571	0.00103470197639755\\
	300	0.00103728205084762\\
	303.571428571429	0.00103988374451228\\
	307.142857142857	0.00104250881879176\\
	310.714285714286	0.00104515910544217\\
	314.285714285714	0.00104783651263781\\
	317.857142857143	0.00105054303348467\\
	321.428571428571	0.001053280754122\\
	325	0.00105605186247168\\
	328.571428571429	0.00105885865841576\\
	332.142857142857	0.0010617035652449\\
	335.714285714286	0.00106458914129539\\
	339.285714285714	0.00106751809348411\\
	342.857142857143	0.00107049329249311\\
	346.428571428571	0.00107351778915732\\
	350	0.00107659483300856\\
	353.571428571429	0.0010797278934394\\
	357.142857142857	0.00108292068255517\\
	360.714285714286	0.00108617718175778\\
	364.285714285714	0.00108950167188012\\
	367.857142857143	0.00109289876604766\\
	371.428571428571	0.00109637344939687\\
	375	0.00109993112178734\\
	378.571428571429	0.00110357764927551\\
	382.142857142857	0.00110731942117383\\
	385.714285714286	0.00111116341811075\\
	389.285714285714	0.00111511728914028\\
	392.857142857143	0.00111918944292514\\
	396.428571428571	0.00112338915455909\\
	400	0.00112772669027599\\
	403.571428571429	0.00113221345691375\\
	407.142857142857	0.00113686217856234\\
	410.714285714286	0.0011416871087019\\
	414.285714285714	0.00114670428692649\\
	417.857142857143	0.00115193184993504\\
	421.428571428571	0.0011573904130968\\
	425	0.00116310354156456\\
	428.571428571429	0.00116909833753818\\
	432.142857142857	0.00117540617876395\\
	435.714285714286	0.00118206365915513\\
	439.285714285714	0.00118911379812231\\
	442.857142857143	0.00119660761863876\\
	446.428571428571	0.00120460623441625\\
	450	0.00121318365512486\\
	453.571428571429	0.00122243062588494\\
	457.142857142857	0.0012324599889361\\
	460.714285714286	0.00124341434603867\\
	464.285714285714	0.00125547729990017\\
	467.857142857143	0.00126889046091361\\
	471.428571428571	0.00128398010601656\\
	475	0.00130120076575913\\
	478.571428571429	0.00132121013499491\\
	482.142857142857	0.00134500561751302\\
	485.714285714286	0.00137418976511899\\
	489.285714285714	0.00141150221534494\\
	492.857142857143	0.00146115465918675\\
	496.428571428571	0.00149603554656411\\
	500	nan\\
};
\addlegendentry{$q$}
\end{axis}
\end{tikzpicture}%
		\caption{``lifted'' controls $q_A$ and $q$  of  \eqref{eq:liftedOCP} using the ZDP with $m=2$. Plotted values are obtained by forward integration with the calculated ``lifted'' controls.}
		\label{fig:cstrZDPControls}
	\emp
\end{figure}
\begin{figure}[ht]
	\centering
	\bmp{0.4}
%\definecolor{mycolor1}{rgb}{0.00000,0.44700,0.74100}%
% all in black
\definecolor{mycolor1}{rgb}{0,0,0}%
\begin{tikzpicture}

\begin{axis}[%
%width = 0.75\textwidth,
width=0.6\textwidth,
at={(0.812in,0.573in)},
scale only axis,
xmin=0,
xmax=500,
xlabel style={font=\color{white!15!black}},
xlabel={time $t$},
xtick={0,100,200,300,400,500},
ymin=0.01,
ymax=0.09,
axis background/.style={fill=white},
axis x line*=bottom,
axis y line*=left,
legend style={at={(0.03,0.97)}, anchor=north west, legend cell align=left, align=left, draw=white!15!black}
]
\addplot [color=mycolor1, line width=1.2pt]
  table[row sep=crcr]{%
0	0.01\\
3.57142857142857	0.0135714565389128\\
7.14285714285714	0.01714287956811\\
10.7142857142857	0.0207142543824583\\
14.2857142857143	0.024285562135626\\
17.8571428571429	0.0278567775967853\\
21.4285714285714	0.0314278655698436\\
25	0.034998774444056\\
28.5714285714286	0.0385694235785611\\
32.1428571428571	0.0421396757739291\\
35.7142857142857	0.0457092671615876\\
39.2857142857143	0.0492775804675063\\
42.8571428571429	0.0528425119303492\\
46.4285714285714	0.0563809438168352\\
50	0.059583346955087\\
53.5714285714286	0.062415925880588\\
57.1428571428571	0.0649390261242773\\
60.7142857142857	0.0671996794903571\\
64.2857142857143	0.0692343030246383\\
67.8571428571429	0.0710716749491993\\
71.4285714285714	0.0727349943545687\\
75	0.0742432903678463\\
78.5714285714286	0.0756124065605354\\
82.1428571428571	0.0768557053873027\\
85.7142857142857	0.0779845837773884\\
89.2857142857143	0.0790088582168751\\
92.8571428571429	0.0799370575723872\\
96.4285714285714	0.0807766493439026\\
100	0.081534216988131\\
103.571428571429	0.0822156006783339\\
107.142857142857	0.0828260103306318\\
110.714285714286	0.0833701173094464\\
114.285714285714	0.0838521295405101\\
117.857142857143	0.0842758535671296\\
121.428571428571	0.0846447462273633\\
125	0.0849619580037912\\
128.571428571429	0.0852303696347413\\
132.142857142857	0.085452623229621\\
135.714285714286	0.0856311488691099\\
139.285714285714	0.0857681874706553\\
142.857142857143	0.085865810545031\\
146.428571428571	0.0859259373495088\\
150	0.085950349848356\\
153.571428571429	0.0859407058168107\\
157.142857142857	0.0858985503647051\\
160.714285714286	0.0858253261085335\\
164.285714285714	0.0857223821818433\\
167.857142857143	0.0855909822427205\\
171.428571428571	0.0854323116117125\\
175	0.0852474836526609\\
178.571428571429	0.0850375454917237\\
182.142857142857	0.0848034831558253\\
185.714285714286	0.0845462261999938\\
189.285714285714	0.0842666518832411\\
192.857142857143	0.0839655889444813\\
196.428571428571	0.0836438210231365\\
200	0.0833020897630845\\
203.571428571429	0.0829410976339879\\
207.142857142857	0.082561510499364\\
210.714285714286	0.0821639599575385\\
214.285714285714	0.0817490454783305\\
217.857142857143	0.0813173363556017\\
221.428571428571	0.080869373493465\\
225	0.0804056710420605\\
228.571428571429	0.0799267178967715\\
232.142857142857	0.079432979073405\\
235.714285714286	0.0789248969703732\\
239.285714285714	0.0784028925277174\\
242.857142857143	0.0778673662917855\\
246.428571428571	0.0773186993933851\\
250	0.0767572544463585\\
253.571428571429	0.0761833763729002\\
257.142857142857	0.0755973931610586\\
260.714285714286	0.0749996165594784\\
264.285714285714	0.0743903427137464\\
267.857142857143	0.0737698527482098\\
271.428571428571	0.073138413296843\\
275	0.0724962769861492\\
278.571428571429	0.0718436828727889\\
282.142857142857	0.0711808568383127\\
285.714285714286	0.0705080119429646\\
289.285714285714	0.0698253487403702\\
292.857142857143	0.069133055554407\\
296.428571428571	0.0684313087195814\\
300	0.0677202727856547\\
303.571428571429	0.0670001006872571\\
307.142857142857	0.0662709338788044\\
310.714285714286	0.0655329024349214\\
314.285714285714	0.0647861251161501\\
317.857142857143	0.0640307093995653\\
321.428571428571	0.0632667514736115\\
325	0.0624943361961187\\
328.571428571429	0.0617135370141549\\
332.142857142857	0.0609244158441202\\
335.714285714286	0.0601270229098798\\
339.285714285714	0.0593213965365414\\
342.857142857143	0.0585075628968001\\
346.428571428571	0.057685535706306\\
350	0.0568553158639682\\
353.571428571429	0.0560168910322355\\
357.142857142857	0.055170235151683\\
360.714285714286	0.054315307883199\\
364.285714285714	0.0534520539701038\\
367.857142857143	0.0525804025110962\\
371.428571428571	0.0517002661333715\\
375	0.0508115400535705\\
378.571428571429	0.049914101012041\\
382.142857142857	0.0490078060631865\\
385.714285714286	0.0480924912018105\\
389.285714285714	0.0471679698013801\\
392.857142857143	0.0462340308358518\\
396.428571428571	0.0452904368510401\\
400	0.0443369216447555\\
403.571428571429	0.0433731876064549\\
407.142857142857	0.0423989026570195\\
410.714285714286	0.0414136967158792\\
414.285714285714	0.0404171576065664\\
417.857142857143	0.0394088262906579\\
421.428571428571	0.0383881912935201\\
425	0.0373546821507279\\
428.571428571429	0.0363076616589563\\
432.142857142857	0.0352464166562633\\
435.714285714286	0.0341701469777989\\
439.285714285714	0.0330779521275672\\
442.857142857143	0.0319688150632534\\
446.428571428571	0.0308415822934915\\
450	0.0296949392110906\\
453.571428571429	0.0285273791960339\\
457.142857142857	0.0273371644623611\\
460.714285714286	0.0261222758088032\\
464.285714285714	0.0248803472331198\\
467.857142857143	0.0236085795847507\\
471.428571428571	0.0223036247785201\\
475	0.0209614282613478\\
478.571428571429	0.0195770124965305\\
482.142857142857	0.0181441809269617\\
485.714285714286	0.0166551373454499\\
489.285714285714	0.0151001514981694\\
492.857142857143	0.0134683567502879\\
496.428571428571	0.011757440119784\\
500	0.01\\
};
\addlegendentry{$V$}
\end{axis}
\end{tikzpicture}%
	 	\caption{``full'' volume $V$ of  \eqref{eq:origOCP} solved by multiple shooting.}
	 	\label{fig:cstrFullVolume}
	\emp
	\qquad 
	\bmp{0.4}
		\input{cstrZDPN140Volume}
		\caption{``lifted'' volume $V$  of  \eqref{eq:liftedOCP} using the ZDP with $m=2$}
		\label{fig:cstrZDPVolume}
	\emp
\end{figure}

The solutions of the original system \eqref{eq:origOCP} via multiple shooting and the lifted system \eqref{eq:liftedOCP} (without forward integration) are, in contrast to the enzyme example, not almost equal. It holds 
%	\begin{align}
%	& \max\left\lbrace	\lvert \lvert 	c^{\text{orig}}-	c^{\text{app}} \rvert \rvert_{\infty},  
%	\lvert \lvert 	q_{\text{A}}^{\text{orig}}-	q_{\text{A}}^{\text{app}} \rvert \rvert_{\infty}, 					
%	\lvert \lvert	q^{\text{orig}}-q^{\text{app}} \rvert \rvert_{\infty} \right\rbrace \nonumber \\
%	&= \lvert \lvert 	c_{\text{B}}^{\text{orig}}-c_{\text{B}}^{\text{app}} \rvert \rvert_{\infty} \approx 0.049.
%	\end{align} 
\[ \max\left\lbrace	\lvert \lvert 	c^{\text{orig}}-	c^{\text{app}} \rvert \rvert_{\infty},  
\lvert \lvert 	q_{\text{A}}^{\text{orig}}-	q_{\text{A}}^{\text{app}} \rvert \rvert_{\infty}, 					
\lvert \lvert	q^{\text{orig}}-q^{\text{app}} \rvert \rvert_{\infty} \right\rbrace 
= \lvert \lvert 	c_{\text{B}}^{\text{orig}}-c_{\text{B}}^{\text{app}} \rvert \rvert_{\infty} \approx 0.049.
\]
This can also be observed in Figures \ref{fig:cstrFullStates} and \ref{fig:cstrZDPStates}. Variable $c_\text{B}$ in the lifted OCP grows more rapidly at the beginning of the time interval. 
IPOPT has found two different optimal solutions (for two different nonlinear problems). The objective function values are $-0.25667$ for the full system and $-0.25679$ for the lifted system which yields an absolute difference of $1.2 \cdot 10^{-4}$ and a relative difference of about $0.047 \%$.  

Averaged runtimes for all methods mentioned above can be found in Table \ref{tab:runtimeCSTR}. Because of the slowness of the reduced OCP, we did not implement this approach for the CSTR example. The relations of the runtimes are comparable to those of the first example. For the lifted systems including the ZDP method of order $m=2$ a significant speed up of factor 38 in comparison with the full multiple shooting method is observed. 

%%% PLACEHOLDER FOR TABLE 3

Depending on the number of intervals, the solution of the full system changes. For $N=4000$ intervals, the solution of the full OCP via multiple shooting (see Figures \ref{fig:cstr4000States},\ref{fig:cstr4000Controls}, and \ref{fig:cstr4000Volume}) looks qualitatively more like the lifted solution for $N=140$ than like the solution of the full OCP. 
The objective function value is $-0.2568167$. These facts point out the sensitivity of the problem as well as the benefits of the lifted OCP. 

%This demonstrates that these problems are very sensitive
%and even small changes (like here in the discretization) lead to different local minima. Furthermore, it shows that the lifted solution is not just an 
%arbitrary solution without practical relevance.

\begin{figure}[ht]
	\centering
	\bmp{0.4}
		\input{cstrFullN4000States}
		\caption{``full'' species $c_A,c_B,c_C$ and $c_D$ of  \eqref{eq:origOCP} solved by multiple shooting on $N=4000$ intervals.}
		\label{fig:cstr4000States}
	\emp
	\qquad
	\bmp{0.4}
		\input{cstrFullN4000Controls}
		\caption{``full'' controls $q_A$ and $q$ of  \eqref{eq:origOCP} solved by multiple shooting on $N=4000$ intervals.}
		\label{fig:cstr4000Controls}
	\emp
\end{figure}
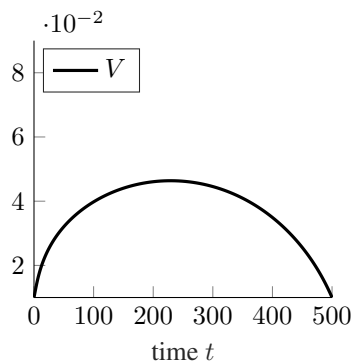
\begin{figure}[ht]
	\centering
	\bmp{0.4}
%\definecolor{mycolor1}{rgb}{0.00000,0.44700,0.74100}%
% all in black
\definecolor{mycolor1}{rgb}{0,0,0}%
\begin{tikzpicture}

\begin{axis}[%
width = 0.6\textwidth,
at={(0.812in,0.573in)},
scale only axis,
xmin=0,
xmax=500,
xlabel style={font=\color{white!15!black}},
xlabel={time $t$},
xtick={ 0, 100, 200, 300, 400, 500},
ymin=0.01,
ymax=0.09,
axis background/.style={fill=white},
axis x line*=bottom,
axis y line*=left,
legend style={at={(0.03,0.97)}, anchor=north west, legend cell align=left, align=left, draw=white!15!black}
]
\addplot [color=mycolor1, line width=1.2pt]
table[row sep=crcr]{%
	0	0.00999999999999091\\
	2	0.0119520387022476\\
	3	0.0129050934809811\\
	3.875	0.0137134593051655\\
	4.625	0.0143790228264606\\
	5.375	0.0150150846854444\\
	6.125	0.0156209251215955\\
	6.875	0.0161977784633791\\
	7.625	0.0167476523607775\\
	8.375	0.0172727199436622\\
	9.25	0.0178567225074744\\
	10.125	0.0184128098654242\\
	11	0.0189436413859312\\
	12	0.0195223497042889\\
	13	0.0200741026103515\\
	14	0.0206014802317895\\
	15.125	0.0211684042091065\\
	16.25	0.0217100165671695\\
	17.5	0.0222849598819153\\
	18.75	0.0228342972482096\\
	20.125	0.0234118064828408\\
	21.5	0.0239638900688419\\
	23	0.0245398824074528\\
	24.5	0.0250909606813821\\
	26.125	0.0256624243975807\\
	27.875	0.0262508606579672\\
	29.625	0.0268138617514637\\
	31.5	0.0273914066976886\\
	33.5	0.027980807406891\\
	35.5	0.0285451405514436\\
	37.625	0.0291196814743557\\
	39.875	0.0297022976245671\\
	42.25	0.0302910532850547\\
	44.75	0.0308841941235869\\
	47.375	0.0314801300631302\\
	50	0.0320508349837496\\
	52.75	0.0326238419355036\\
	55.625	0.03319781639334\\
	58.625	0.033771525090458\\
	61.75	0.0343438222758436\\
	65	0.0349136366796188\\
	68.375	0.0354799593597477\\
	71.875	0.0360418325325895\\
	75.5	0.0365983394273144\\
	79.25	0.0371485951633304\\
	83.125	0.0376917386159903\\
	87.125	0.0382269252132232\\
	91.25	0.0387533205907289\\
	95.5	0.0392700950231983\\
	99.875	0.0397764185441929\\
	104.375	0.0402714566656641\\
	109	0.0407543666078709\\
	113.75	0.0412242939536327\\
	118.625	0.0416803696426769\\
	123.5	0.0421109913365854\\
	128.5	0.0425272458422796\\
	133.625	0.042928181115883\\
	138.875	0.0433128247074137\\
	144.125	0.0436719356255821\\
	149.5	0.0440139414699843\\
	155	0.0443377836699597\\
	160.5	0.0446358946531404\\
	166.125	0.0449148047925974\\
	171.75	0.0451680232155809\\
	177.5	0.0454008478748165\\
	183.25	0.0456078435366294\\
	189	0.0457894312677922\\
	194.875	0.0459491043048956\\
	200.75	0.0460829696676228\\
	206.625	0.0461913114191361\\
	212.5	0.0462743657906231\\
	218.5	0.0463332856818397\\
	224.5	0.0463661878183643\\
	230.5	0.0463731857219045\\
	236.5	0.0463543502756352\\
	242.5	0.0463097109943647\\
	248.5	0.0462392569423287\\
	254.5	0.0461429373265219\\
	260.5	0.0460206617866561\\
	266.5	0.0458723003960131\\
	272.375	0.0457015903770639\\
	278.25	0.0455055105385327\\
	284.125	0.0452838260623594\\
	290	0.0450362632117844\\
	295.75	0.0447686077038725\\
	301.5	0.0444755296656467\\
	307.25	0.0441566559888997\\
	312.875	0.0438193552098483\\
	318.5	0.0434565412895154\\
	324.125	0.0430677414431102\\
	329.625	0.0426619579350245\\
	335.125	0.0422303047525929\\
	340.5	0.0417829131865233\\
	345.875	0.0413096746524388\\
	351.125	0.0408218776270814\\
	356.375	0.0403081658924975\\
	361.5	0.0397810176428379\\
	366.5	0.0392416064718191\\
	371.5	0.0386766668628979\\
	376.375	0.0381005343708125\\
	381.25	0.0374986233230743\\
	386	0.0368865475242046\\
	390.625	0.0362654897576817\\
	395.25	0.0356188225906067\\
	399.75	0.0349641769511209\\
	404.125	0.0343027615240317\\
	408.5	0.0336158124382564\\
	412.75	0.0329230880399791\\
	416.875	0.0322258370418353\\
	421	0.0315030440300461\\
	425	0.0307767246814592\\
	428.875	0.0300481828672901\\
	432.75	0.0292940140058136\\
	436.5	0.0285386445735298\\
	440.125	0.0277834511288688\\
	443.75	0.0270024772457305\\
	447.25	0.0262227425860715\\
	450.625	0.0254457184382773\\
	454	0.0246427065316652\\
	457.25	0.0238435384550826\\
	460.375	0.0230498104364187\\
	463.375	0.0222631849981099\\
	466.375	0.0214510267146579\\
	469.25	0.0206473432378402\\
	472	0.0198539867535601\\
	474.75	0.0190351054428675\\
	477.375	0.018228181741506\\
	480	0.0173950458155332\\
	482.5	0.016575704245156\\
	484.875	0.0157725202781194\\
	487.25	0.0149437388657248\\
	489.5	0.0141336517069135\\
	491.75	0.0132980598127119\\
	494	0.0124358099843676\\
	496.125	0.01159611226268\\
	498.25	0.010731240462917\\
	500	0.00999999999999091\\
};
\addlegendentry{$V$}
\end{axis}
\end{tikzpicture}%
		\caption{``full'' volume $V$ of  \eqref{eq:origOCP} solved by multiple shooting on $N=4000$ intervals.}
		\label{fig:cstr4000Volume}
	\emp
\end{figure}

\section{Conclusion}
\label{sec:conclusion}
In this paper, we show by help of benchmark problems that it is possible to use slow manifold based model reduction techniques in order to 
solve optimal control problems (OCPs) with multiple time-scale dynamics based on multiple shooting faster but slightly less accurate. 
For this reason, the dynamics of the fast optimization variables are replaced by nonlinear equality conditions.
The benefit of that approach is not visible at first sight since the system we get in this way is of the same 
order -- so we cannot speak of model order reduction in proper meaning of the word. However, the remaining 
ODE part of the OCP is much less stiff. This allows us to use explicit integrating routines rather
than implicit integrators to solve the ODE part numerically.

Numerical experiments have shown that this approach gives a speed-up of factor 7 to 38. Since the new approach differs from the original OCP, in general different solutions are expected. But, the corresponding objective function values differ just slightly. 

%	In comparison with direct collocation, our proposed approach is a little bit slower, but uses significantly less optimization variables. 
%	Probably, much higher-dimensional OCP models could significantly benefit from our proposed method, since direct collocation leads to very large problems.

A goal for future research is to improve the new approach in order to be able to solve OCPs with stiff chemical kinetic ODE constraints in real-time. 

% set path of bib-file correctly
\newpage
\bibliography{literature}

\newpage
\appendix

\section{Tables}

\begin{table}[htp]
	\caption{Runtimes in seconds, the number of variables and constraints for the underlying NLP of the enzyme example with $N=40$.}
	\label{tab:runtimeEnzyme}
	\centering
	\begin{tabular}{lrrr}
		{\bfseries method} 											& {\bfseries runtime} 	& {\bfseries \#variables} & {\bfseries \#constraints}\\  \hline 
		\cellcolor{cell}full OCP (multiple Shooting) 		& \cellcolor{cell}0.1688 & \cellcolor{cell}120  	& \cellcolor{cell}80 \\
%			full OCP (direct Collocation) 							&  0.0313					& 280 							& 240 \\
		reduced OCP	(ZDP $m=2$)							& 0.5272 						& 80 							& 40\\
		\cellcolor{cell}lifted OCP (Lebiedz/Unger)		& \cellcolor{cell} 0.0268	 & \cellcolor{cell} 121		&\cellcolor{cell} 80\\
		lifted OCP (ZDP $m=2$)								&	0.0224 						& 121 						& 80
	\end{tabular}
\end{table}

\begin{table}[ht]
	\caption{constants for the CSTR example}
	\label{table:cstrConstants}
	\centering
	\begin{tabular}{lll}
		\multicolumn{1}{c}{\bfseries{constant}}	& \multicolumn{1}{c}{\bfseries{description}}
		& \multicolumn{1}{c}{{\bfseries value}} \\ \hline
		$c_{A_{\text{in}}}$ & concentration of species $A$ of inlet 						& $1$ \\
		$k_1$				& forward reaction rate constant of $A \rightleftharpoons B$ 			& $100$ \\
		$k_1^-$				& backward reaction rate constant of $A \rightleftharpoons B$			& $90$ \\
		$k_2$				& reaction rate constant of $C \rightarrow D$ 							& $20$ \\
		$\epsilon$			& reaction rate constant of $B \rightarrow C$ 							& $10^{-6}$	\\	
		$q_{A_\text{max}}$	& maximum value for $q_{A}$ in inlet 						& $10^{-3}$\\ 
		$q_{\text{max}}$	& maximum value for $q$ in outlet						& $1.5 \cdot 10^{-3}$\\
		$c_A^*$				& initial concentration of species $A$							& $10^{-3}$\\
		$c_B^*$				& initial concentration of species $B$							& $10^{-3}$\\
		$c_C^*$				& initial concentration of species $C$							& $0$\\
		$c_D^*$				& initial concentration of species $D$							& $10^{-8}$\\
		$V^*$				& initial and final volume of tank reactor						& $10^{-2}$\\
	\end{tabular}
\end{table}

\begin{table}[htp]
	\caption{Runtimes for the CSTR example in seconds for $N=140$.}
	\label{tab:runtimeCSTR}
	\centering
	\begin{tabular}{lrrr}
		{\bfseries method} 											& {\bfseries runtime} 		& {\bfseries \#variables} & {\bfseries \#constraints}\\  \hline 
		\cellcolor{cell}full OCP (multiple Shooting) 		& \cellcolor{cell}5.0817 	& \cellcolor{cell}979 		 & \cellcolor{cell}700 \\
%			full OCP (direct Collocation) 							&  0.1312						& 2379 							& 2100 \\
		lifted OCP (Lebiedz/Unger)								&0.2473	 						&980							 &700\\
		\cellcolor{cell}lifted OCP (ZDP $m=2$)				&\cellcolor{cell} 0.1311	& \cellcolor{cell}980  		&\cellcolor{cell} 700
	\end{tabular}
\end{table}

\end{document}